\documentclass[11pt]{article}
\usepackage{amsmath}
\usepackage{amssymb}
\usepackage{dsfont} % for unit operator  (indicator function) \mathds{1}

\usepackage{graphicx}
\usepackage[dvipsnames]{xcolor}
\newcommand{\be}{\begin{equation}}
\newcommand{\ee}{\end{equation}}
\newcommand{\bel}[1]{\begin{equation}\label{#1}}
\newcommand{\bea}{\begin{eqnarray}}
\newcommand{\eea}{\end{eqnarray}}
\newcommand{\balign}{\begin{align}}
\newcommand{\ealign}{\end{align}}
\newcommand{\ba}{\begin{array}}
\newcommand{\ea}{\end{array}}
\newcommand{\bfig}{\begin{figure}}
\newcommand{\efig}{\end{figure}}

\newcommand{\eref}[1]{(\ref{#1})}

\allowdisplaybreaks

\newcommand{\bra}[1]{\mbox{$\langle \, {#1}\, |$}}
\newcommand{\ket}[1]{\mbox{$| \, {#1}\, \rangle$}}

\newcommand{\inprod}[2]{\mbox{$\langle \, {#1} \, | \, {#2} \, \rangle$}}
\newcommand{\Prob}[1]{\mbox{${\rm Prob}\left[ \, {#1}\, \right]$}}

\newcommand{\E}{{\mathbf E}}

\newcommand{\bfx}{\mathbf{x}}
\newcommand{\bfy}{\mathbf{y}}

\newcommand{\floor}[1]{\lfloor{#1}\rfloor}
\newcommand{\ceil}[1]{\lceil{#1}\rceil}

\newcommand{\rmd}{\mathrm{d}}
\newcommand{\rme}{\mathrm{e}}

\newcommand{\sgn}{\mathrm{sgn}}

\newcommand{\ddt}{\frac{\rmd}{\rmd t}}

\newcommand{\comm}[2]{\mbox{$[\,{#1}\,,\,{#2}\,]$}}

\newcommand{\C}{{\mathbb C}}
\newcommand{\R}{{\mathbb R}}

\newcommand{\Z}{{\mathbb Z}}
\newcommand{\N}{{\mathbb N}}
\renewcommand{\S}{\mathbb S}
\newcommand{\bzeta}{\boldsymbol{\zeta}}
\newcommand{\bfeta}{\boldsymbol{\eta}}

\newcommand{\bxi}{\boldsymbol{\xi}}

\newtheorem{theo}{Theorem}[section]
\newtheorem{lmm}[theo]{Lemma}
\newtheorem{df}[theo]{Definition}
\newtheorem{prop}[theo]{Proposition}
\newtheorem{cor}[theo]{Corollary}
\newtheorem{rem}[theo]{Remark}
\newcommand{\Proof}{\noindent {\it Proof: }}
\def\qed{\hfill$\Box$\par\medskip\par\relax}

\begin{document}

\title{Self-duality and shock dynamics in the 
$n$-component priority ASEP}
\author{V. Belitsky$^{1}$ 
%\and Author B$^{1}$
\and G.M.~Sch\"utz$^{2,3}$ 
}

\maketitle
%\today

{\small%\footnotesize
\noindent $^{~1}$Instituto de Matem\'atica e Est\'atistica,
Universidade de S\~ao Paulo, Rua do Mat\~ao, 1010, CEP 05508-090,
S\~ao Paulo - SP, Brazil
\\
\noindent Email: belitsky@ime.usp.br  
%seguchi@ime.usp.br
%\noindent url: \texttt{http://www.proba.jussieu.fr/$\sim$comets}

\smallskip
\noindent $^{~2}$Institute of Complex Systems II,
Forschungszentrum J\"ulich, 52425 J\"ulich, Germany
\\
\noindent Email: g.schuetz@fz-juelich.de

\smallskip
\noindent $^{~3}$Interdisziplin\"ares Zentrum f\"ur Komplexe Systeme, Universit\"at
Bonn, Br\"uhler Str. 7, 53119 Bonn, Germany
}

\begin{abstract}
We study the $n$-component priority asymmetric simple exclusion process 
($n$-ASEP)
with reflecting boundaries. We obtain all invariant measures in explicit form and
prove reversibility. Using the symmetry of the generator of the process under
the quantum algebra $U_q[\mathfrak{gl}(n+1)]$ we construct duality functions with 
respect to which the $n$-ASEP is self-dual, both for the finite and the infinite
integer lattice. For the $n$-ASEP on the infinite lattice we use self-duality to 
derive in explicit form the time evolution of a family of measures with $K$ 
shocks in terms of the transition probability of $K$ coloured particles in a shock 
exclusion process with particle-dependent hopping rates and nearest-neighbour 
colour exchange. This process is a gas of particles that forms a bound state,
corresponding to shock coalescence on macroscopic scale.
\\[.3cm]\textbf{Keywords:} Asymmetric simple exclusion process, Duality,
Quantum algebras, Shocks
\\[.3cm]\textbf{AMS 2000 subject classifications:} 60K35. Secondary: 82C20, 82C23
%\\[.3cm]\textbf{PACS numbers:} 

\end{abstract}

\newpage

\section{Introduction}

\subsection{Informal overview: Process, main tools, and results}

We consider an asymmetric simple exclusion process (ASEP) with $n$ 
species of particles and a priority jump rule where particles
of a species $\alpha\in \S_{0,n} := \{0,\dots,n\}$ 
``see'' particles of a lower species $\beta<\alpha$
as vacant sites \cite{Alca93,Alca00,Arit09,Trac13}.
The Markov dynamics of this process local occupation variables 
$\eta_k \in \S_{0,n}$ on a finite one-dimensional integer lattice 
$\Lambda^{\pm} = [L^-,L^+] \cap \Z$ with 
$L=L^+-L^-+1$ sites can be described informally as follows. Each site 
$k\in \Lambda^{\pm}$ can be either empty (denoted by $0$) or occupied by at 
most particle of species $\alpha$ with $1 \leq \alpha \leq n$. 
Each bond $(k,k+1)$ of $\Lambda^{\pm}$, $L^-\leq k \leq L^+-1$
carries a clock $k$. If $\eta_k\neq\eta_{k+1}$ the clock 
rings independently of all other clocks 
after an exponentially distributed
random time with parameter $\omega_k$ where $\omega_k=wq>0$ if $\eta_k>\eta_{k+1}$ 
and $\omega_k=wq^{-1}>0$ if $\eta_k<\eta_{k+1}$. 
When the clock rings the particle occupation variables are
interchanged and the clocks $k-1$, $k$ and $k+1$ instantly acquire 
the corresponding new parameter.
If $\eta_k=\eta_{k+1}$ then nothing happens, corresponding to parameter $\omega_k=0$.
We shall refer to this process as $n$-priority ASEP.
We consider (i) the finite system with reflecting boundary conditions, 
which means that no jumps from the left boundary site $L^-$ to the left and 
no jumps from the right boundary site $L^+$ to the right are allowed, (ii)
the semi-infinite system with $L^+\to\infty$ or $L^-\to-\infty$, and (iii)
the infinite system defined on $\Z$.

Besides standard probabilistic tools for stochastic interacting 
particle systems a convenient method to discuss this
process is the so-called quantum Hamiltonian approach
to interacting particle systems, described in probabilistic terms
by Sudbury et al. in \cite{Sudb95,Lloy96}. This method reveals \cite{Alca93}
that the generator of the $n$-species priority ASEP is related
to the Hamiltonian operator of an integrable quantum spin system, 
viz. the Perk-Schultz chain 
\cite{Perk81}, which is symmetric under the 
action of the quantum algebra $U_q[\mathfrak{gl}(n+1)]$, 
i.e., the $q$-deformed universal enveloping algebra of 
the Lie algebra $\mathfrak{gl}(n+1)$ \cite{Jimb85,Jimb86}.

The three main results that we obtain from this approach are the following.
(i) {\it Reversibility:} We present in explicit form all invariant 
measures for the finite system (Theorem \eref{Theo:revmeas}). 
We prove reversibility using the detailed balance condition and obtain
the normalization factor from combinatorial arguments.
We also construct some blocking measures for the infinite system.

(ii) {\it Self-duality:} We construct by a ground state transformation the 
representation matrices of $U_q[\mathfrak{gl}(n+1)]$ that commute with the 
generator of the process. Using general probabilistic arguments relating 
reversibility, symmetry and self-duality 
\cite{Schu94,Giar09} we obtain from these self-duality functions 
(Theorem \eref{Theo:duality}). The proof is constructive
and can be generalized to obtain other duality functions.
Also duality functions for the infinite systems, where some
convergence issues need to be taken into account, are derived. 

(iii) {\it Microscopic structure of shocks:} This 
is the main application of duality in this work.
We first describe informally some features of the process in the hydrodynamic 
limit, particularly the appearance of shocks 
which until now have remained elusive in the 
rigorous treatment of stochastic particle systems with more than one 
conservation law due to the lack of attractiveness, see e.g. 
\cite{Toth03,Frit04,Toth05,Popk12}. For the $n$-priority ASEP on $\Z$ 
we then define shock measures with $K$ consecutive shocks of species $n$, 
marked microscopically by $K$ particles of species $\alpha<n$. We prove 
that the time evolution of these shock measures can be expressed in terms 
of the transition probabilities of a different exclusion process on $\Z$ 
with only $K$ particles and particle-dependent 
hopping rates which can be interpreted from a physics perspective as random 
walkers with different masses and on-site repulsion in a constant gravitational 
field. The $K$-particle transition probabilities of this 
``shock exclusion process'' can be computed from the nested Bethe ansatz 
\cite{Yang67,Arit09,Trac13,Gaud14}. 
A single shock, in particular, performs a biased random walk.
The stationary microscopic distances between the shock markers
are independent geometrically
distributed random variables, thus elucidating the microscopic meaning of {\it macroscopic}
coalescence of shocks \cite{Ferr00}. Moreover, the result exposes a link between
the time evolution of shocks, current fluctuations in the ASEP 
\cite{Imam11,Cari16}, and bound states in quantum spin systems \cite{Gaud14}.

\subsection{Setting of the problem}

This work makes use of connections between probability theory,
non-equilibrium statistical mechanics and integrability that have been known 
for a long time, but which have been explored more intensely only recently. 
In order to place our results relating the concepts of duality, symmetry and 
shocks into this context we first mention that
the $n$-priority ASEP is a natural generalization of the standard ASEP ($n=1$) 
\cite{Ligg85,Ligg99,Schu01} to several conserved species of particles.
The 2-priority ASEP is the ASEP with second-class particles going back to
\cite{Ligg76} and studied recently in the context of duality in
\cite{Beli15a,Beli15b,Kuan16} for finite lattices with reflecting
boundary conditions. During the final stages of this work
we were notified of a closely related result where a duality 
function of an $n$-priority ASEP with up to
$2j$ particles per site is derived \cite{Kuan16b}. For $j=1/2$ this duality 
function is related to the duality function of Theorem \eref{Theo:duality}
but slightly different, see Remark \eref{Rem:duality}.

It should be noted that reflecting boundary conditions lead to rather 
different properties of the process than periodic boundary conditions where 
the invariant measure, which can be expressed in a matrix product form
in the totally asymmetric limit $q\to\infty$ \cite{Mall99,Ferr07,Evan09,Prol09},
is not reversible for $q\neq 1$. Moreover, the $U_q[\mathfrak{gl}(n+1)]$-symmetry 
is broken for $q\neq 1$, except possibly for an unexplored residual property 
known for the simplest case $n=1$ \cite{Pasq90,Schu16}. The invariant measures 
that we obtain for reflecting boundary conditions
turn out to have the peculiar property of long-range interactions despite
local nearest-neighbour dynamics, reminiscent of a similar phenomenon found
in another exclusion process with periodic boundary conditions
\cite{Evan98b,Clin03,Ayer09,Bodi11}. An intriguing question is whether
the matrix product measures for the periodic system are related to
(and can perhaps be constructed from) some residual quantum algebra symmetry.
Also self-duality for open boundaries where first-class particles are injected 
and extracted but second-class particles are reflected
\cite{Kreb03,Cram15,Cram16,Mand16} is an open problem.

The idea of using non-Abelian symmetries of the generator for deriving duality 
relations for stochastic interacting particle systems goes back to 
\cite{Schu94,Schu97} where it was shown 
that duality functions arise from representations 
of the symmetry algebra. The range of models that can be treated in 
this fashion is large since non-Abelian symmetries, in particular 
Lie algebras and their quantum deformations \cite{Jimb86}, appear frequently 
in integrable quantum systems and some of their non-integrable generalizations, 
many of which are related to generators of stochastic interacting particle 
systems \cite{Alca93,Schu01,Giar09}. Thus, given a symmetry, the derivation 
of duality functions reduces to finding those representations of the symmetry
algebra that commute with the generator of the stochastic 
interacting particle system and to computing the matrix elements
of representations of the symmetry operator.
This approach was brought into a neat and 
systematic form by Giardin\`a et al. \cite{Giar09} and Jansen and Kurt
\cite{Jans14} and was applied to various interacting particle systems
to study current fluctuations, shock motion, heat conduction 
and other properties of these systems 
\cite{Giar07,Ohku10,Imam11,Cari13,Boro14,Cari15,Corw15,Cari16}.

The two-species priority ASEP on $\Z$ has been used to study 
microscopic properties of shocks \cite{Ferr91,Derr93,Ferr94a,Ferr94b}. 
Particles of species 2 are then the so-called first-class particles and particles 
of species 1 (called second-class
particles) serve to define the microscopic positions of shocks which 
become macroscopic density discontinuities in the large-scale
behaviour of usual ASEP ($n=1$). An intriguing property
is the fact that the time-evolution of certain shock measures
with $k$ shocks for configurations with an arbitrary number of particles
is given by the transition probability of $K$ exclusion
particles \cite{Beli02,Bala10}. This fact, which appears to be related
to the existence of bound states in the associated quantum system, 
provides detailed information
microscopic properties of these shocks. The present work demonstrates  that 
the $K$-particle property of these shock measures 
for the case of $n=1$ and $n=2$ arise from self-duality and 
shows that a similar property holds for general $n$. Conversely, random walk
properties of shocks have been proved for various other processes 
\cite{Paes04,Rako04,Bala10}. This may allow for
finding dualities and non-Abelian symmetries in these processes.
Our results also indicate a link between current fluctuations \cite{Sasa10,Prol10,Amir11,Imam11,Laza11,Gori12,Cari16} 
and the dynamics of shocks
\cite{Pigo00,Beli02,Bala10,Beli11,Civi14} via self-duality
since for both problems the same duality functions are used.

\subsection{Structure of the paper}

In Sec. 2 we define the generator of the
process and collect some known facts required to state and prove the main results, 
which are presented in Sec. 3 along with some further remarks on properties 
of the $n$-component priority ASEP. In Sec. 4 we present the proofs.

\section{The $n$-component priority ASEP}

\subsection{Definitions and notation}
\label{Sec:Definotat}

Here we fix notation necessary to define the $n$-priority ASEP. 
Various other conventions and frequently used formulas
are collected in Appendix \eref{A1}.

\subsubsection{State space and configurations}
\label{Sec:Defconfig}

It is expedient to employ two distinct
representations of the configurations of the the $n$-priority ASEP that
we shall call the {\it occupation variable representation} $\bfeta$, 
which specifies the type of particle $\eta(k)$ on any given lattice site $k$, 
and the {\it coordinate representation} which specifies the positions $x_i$
and species $\alpha_i$ of the particles which are all uniquely tagged with a 
label $i$. To make these representations and their relation with each other
precise we introduce first
the lattices and the sets of particle species
that we shall work with.

\begin{df}
\label{Def:lattice}
For integers $L^\pm$ with $L:=L^+ +1 -L^- \geq 2$ 
we define four integer lattices $\Lambda^{\cdot}$ by the sets
\bel{lattices} 
\Z; \quad
\Lambda_+ := \{ k\in \Z: k \geq L^-\} ; \quad
\Lambda_- := \{ k\in \Z: k \leq L^+\}; \quad 
\Lambda_{L} := \Lambda_- \cap \Lambda_+ .
\ee
An element $k \in \Lambda$ is referred to as a site with coordinate $k$.
Elements of the sets of integers
\bel{localset}
\S_{n} := \{1,\dots,n\}, \quad \S_{0,n} := \{0,1,\dots,n\}
\ee
are called particle species.
\end{df}

\noindent \underline{\it A) Occupation variable representation:}\\

For a given lattice $\Lambda$ we denote by $\bfeta = (\eta_{k})_{k\in \Lambda}$ 
the configuration of the particle system where lattice
site $k\in\Lambda$ is occupied by a particle of species 
$\eta_k \in \S_{0,n}$. The set of configurations is denoted
$\S_{0,n}^{\Lambda}$ with the short-hands $\Lambda \to L$ for $\Lambda = 
\Lambda_{L}$ and $\Lambda\to\pm$ for $\Lambda = \Lambda_{\pm}$. 
We refer to particles of species 
$0$ also as vacancies.

Next we define the following subsets of $\S_{0,n}^{\Lambda}$.

\begin{df}
\label{Def:subsets}
Let $\Lambda$ be one of the lattices defined in \eref{Def:lattice} and 
define for $\bfeta\in \S_{0,n}^{\Lambda}$ the particle numbers
\be 
\label{Def:partnum}
N^\alpha(\bfeta) := \sum_{k\in\Lambda} \delta_{\eta_k,\alpha}, \quad
M^\alpha(\bfeta) := \sum_{\beta=\alpha}^n N^\beta(\bfeta), \quad \alpha\in \S_{0,n}
\ee
and the projectors 
\bel{particleproj}
\wp^\alpha_{\vec{N}} (\bfeta) := 
\prod_{\beta=\alpha}^n \delta_{N^\alpha(\bfeta),N^\alpha}
\ee
on particle numbers $N^\beta$ with $\beta \geq \alpha$
and denote by $N(\bfeta):= M^1(\bfeta)$ the total particle number
and $\mathcal{P}_{\vec{N}}(\bfeta) := \wp^1_{\vec{N}}(\bfeta)$
the projector on particle numbers $N^1,\dots,N^n$. \\

\noindent (a) For specified particle number $N^\alpha \geq 0$ the subset
\bel{finiteNa}
\S_{N^\alpha}^\Lambda := 
\{\bfeta \in \S_{0,n}^\Lambda \, | \, N^\alpha(\bfeta) = N^\alpha \}
\ee
of $\S_{0,n}^\Lambda$ is called the set of configurations with particle 
number $N^\alpha$ of species $\alpha$ and
for specified particle numbers $\vec{N} = (N^1,\dots,N^n)$ 
with $N^\alpha\geq 0$ the subset
\bel{finiteN}
\S_{\vec{N}}^\Lambda := 
\{\bfeta \in \S_{0,n}^\Lambda \, | \, \mathcal{P}_{\vec{N}} (\bfeta)=1\}
\ee
of $\S_{0,n}^\Lambda$ is called the set of configurations with particle 
numbers $N^\alpha$ of species $\alpha \geq 1$.\\

\noindent (b) The sets of configurations which are right- or left-asymptotically 
fully occupied by particle of species $\alpha\in\S_{0,n}$ are defined by
\bea  
\label{rfull}
\S^{\Lambda}_{\alpha^{>}} & := &
\{\bfeta\in \S_{0,n}^{\Z}: 
\sum_{k=1}^\infty (1-\delta_{\eta_k,\alpha}) < \infty\} \\ 
\label{lfull}
\S^{\Lambda}_{\alpha^{<}} & := &
\{\bfeta\in \S_{0,n}^{\Z}: 
\sum_{k=-\infty}^0 (1-\delta_{\eta_k,\alpha}) < \infty\}.
\eea
\end{df}
For the finite lattice $\Lambda^L$ one has trivially $\S^{L}_{\alpha^{>}} = \S^{L}_{\alpha^{<}}
= \S_{0,n}^L$ and for the semi-infinite lattices $\Lambda^\pm$ one has 
$\S^{+}_{\alpha^{<}} = \S_{0,n}^+$ and
$\S^{-}_{\alpha^{>}} = \S_{0,n}^-$ for any $\alpha$.\\

\noindent \underline{\it B) Coordinate representation:}\\

Configurations in $\S^{\Lambda}_{0^{<}}$ which are left-asymptotically
vacant
%(and hence all configurations on $\Lambda^L$ and $\Lambda^+$) 
can be specified in an alternative way by
consecutively indexing each particle of species $\alpha \geq 1$ with an integer
$i$ such that the left-most particle (which may be of any species $\alpha \geq 1$) 
is assigned the index $1$. The sites occupied by particles
of species $\alpha\geq 1$ are denoted $x^\alpha_i$ with $1\leq i \leq N^\alpha$
and $\{x^\alpha\}$ is the set of these sites. The set 
of all sites occupied by a particle of species $\alpha\geq 1$ is
$\{x_1,\dots,x_N\} = \cup_{\alpha>0} \{x^\alpha\}$. The colour
of particle $i$ is denoted by $\alpha_i$ with $1\leq i \leq N$.

Many applications,
such as the Bethe ansatz, are based on this coordinate representation,
which is also frequently used below. 
The $n$-species priority ASEP becomes in this language an
exclusion process which can be informally described as follows: 
Particles carry a ``colour'' index $\alpha \in \S_n$,
corresponding to species $\alpha \geq 1$ in the language of the
$n$-priority ASEP. If site $x_i\pm 1$ is empty, particle $i$ jumps
after an exponentially distributed random time with parameter
$w q^{\pm 1}$.
If two particles $i,i+1$ are nearest neighbours they do not jump but
exchange their colour after an exponential random time 
with parameter $q$ if 
$n \geq \alpha_j > \alpha_{j+1}\geq 1$, and with parameter 
$q^{-1}$ if $1\leq \alpha_j < \alpha_{j+1} \leq n$. 
Notice that the order $x_{j+1} > x_j$ of the particle coordinates
remains preserved. 

To define the state space in the coordinate representation formally we 
introduce some more notions.

\begin{df}
\label{Def:Weylalcove} 
\noindent a) (Weyl alcove of type $\tilde{C}_{N}$ \cite{Grab99,Grab02,Krat07}) 
Let $\Lambda$ be one of the lattices defined in \eref{Def:lattice}.
For a strictly positive integer $N$ we define the (shifted and scaled) 
Weyl alcove $W^\Lambda_N$  by the set of coordinate vectors 
$\vec{x} := (x_1,\dots,x_N) \in \Lambda^N$ satisfying
$L^- \leq x_1 < x_2 < \dots < x_{N} \leq L^+$.
For $N=0$ we define
$W^\Lambda_0 := \emptyset$. We also define 
\be 
\label{fullcoordinateset}
W^{\Lambda} = \underset{N \geq 0}{\bigcup} W^{\Lambda}_{N}.
\ee
and 
for any $N\geq 1$ the coordinate set 
$\{ \vec{x} \} = \{x_1,\dots,x_N\} \subset \Lambda$  \\

\noindent b) (Colour array) For $\alpha_i \in \S_n$ we define the $N$-particle
colour array as $\vec{\alpha} := (\alpha_1,\dots,\alpha_N) \in \S_n^N$.
We define also the sets 
\be 
\label{fullconfigset}
V^\Lambda_N := W^\Lambda_N \times \S_n^N, \quad
V^{\Lambda} = \underset{N \geq 0}{\bigcup} V^{\Lambda}_{N}
\ee
and denote elements of $V^\Lambda_N$ by the pair $\bfx = (\vec{x},\vec{\alpha})$.
\end{df}

For a specific lattice $\Lambda$ the corresponding Weyl alcoves are denoted 
by $W^L_N$, $W^\pm_N$, and $W^{\infty}_N$ respectively and similarly for 
$V^{\Lambda}$.\\

\noindent \underline{\it Relation between occupation variable and coordinate representations:}\\

We have an isomorphism between configurations
$\bfx \in V^\Lambda_N$ and $\bfeta \in \S^{\Lambda}_{\vec{N}}$
through the bijection $\bfx \leftrightarrow \bfeta$ where 
\be 
\label{isoconfig}
\eta_{k}(\bfx) = \sum_{i=1}^{N(\bfx)}
\alpha_i \delta_{k,x_i}, \quad
x_i (\bfeta) 
= \min_k \left(\sum_{l=L^-}^k (1-\delta_{\eta_l,0}) = i\right), \,
\alpha_i (\bfeta) = \eta_{x_i}
\ee
Defining $\mathbf{0}=(\dots,0,0,0,\dots) \in \S^{\Lambda}_{\alpha^{<}}$
as representing configuration corresponding to the empty lattice this bijection 
together with the bijection $\bfx(\mathbf{0}) = \emptyset$ and
$\bfeta(\emptyset) = \mathbf{0}$ 
yields also an isomorphism between $\bfx \in W^\Lambda \times \S_n^\N$ and 
$\bfeta \in \S^{\Lambda}_{0^{<}}$ for infinite $N$. 

In some instances it will indeed be convenient to view the particle positions as 
a function of a configuration $\bfeta$ and vice versa. 
The function $\bfx: \S_{0^{<}}^\Lambda \mapsto V^{\Lambda}$
is then denoted by $\bfx(\bfeta)$,
Conversely $\bfeta(\bfx)$ is the occupation variable representation interpreted as 
a function $\bfeta: V^\Lambda \mapsto \S_{0^{<}}^\Lambda$ 
of a particle configuration $\bfx$ into the
occupation variable representation $\bfeta(\bfx)$. 
We shall also use
the function $\vec{x}^\alpha: \S_{0,n}^\Lambda \mapsto W^{\Lambda}_{N^\alpha}$
which gives only coordinates $\vec{x}^\alpha(\bfeta)$ 
of particles of species $\alpha$ and the function $x^\alpha_i: 
\S_{0,n}^\Lambda \mapsto \Lambda$ gives the particle position $x^\alpha_i(\bfeta)$ 
of the $i^{th}$ particle of species $\alpha$. 

Unless stated otherwise, $\sum_{\bfeta}$ will always be understood as 
$\sum_{\bfeta \in \S_{0,n}^\Lambda}$ and correspondingly 
$\sum_{\bfx} \equiv \sum_{\bfx\in V^{\Lambda}}$ and similar for products.
A sum over an empty index set is defined to be zero and a product over an
empty index set is defined to be equal to 1.

\subsubsection{Functions of the configurations}

Several other functions of the configurations $\bfeta \in \S^\Lambda_{0,n}$ will 
play a role. Through the bijection
\eref{isoconfig} these induce analogous mappings in the coordinate
representations that we do not all write explicitly. Generally, however,
one has from \eref{isoconfig}
\be 
\delta_{\eta_k,\alpha} = \delta_{\alpha,0} + \sum_{i=1}^{N(\bfx)} \delta_{x_i,k}
(\delta_{\alpha,\alpha_i} - \delta_{\alpha,0})
\ee
and therefore for any function $f(\eta_k)$ of the local occupation variable
\bel{isofun}
f(\eta_k) = f(0) + 
\sum_{i=1}^{N(\bfx)} \delta_{x_i,k} (f(\alpha_i)-f(0))
\ee
where the r.h.s. is expressed in terms of 
the coordinate representation $\bfx(\bfeta)$. By iteration one
obtains form this formula analogous expressions for arbitrary cylinder 
functions.

\begin{df}
\label{Def:functions}
Let $\bfeta \in \S_{0,n}^{\Lambda}$ where $\Lambda$ is one of the
four lattices defined in \eref{Def:lattice}.\\

\noindent a) For $k \in \Lambda$ the local cyclic flip operation 
$\gamma_k(\bfeta)$ 
is defined by
\bel{flip}
(\gamma_k(\bfeta)_l = \left\{ \ba{ll} \eta_k + 1 \ \mathrm{mod} \
(n+1) & l=k \\
\eta_l & l \neq k \ea \right. 
\ee
and $\Gamma := \gamma_{L^-} \circ \dots \circ \gamma_{L^+}$ is called 
cyclic flip. For a given $\bfeta$ we write $\bfeta_k^{p} := 
(\gamma_k)^p(\bfeta)$ and $\bfeta^{p} := \Gamma^p(\bfeta)$
for the $p$-fold action of $\gamma_k$ and $\Gamma$ resp.\\

\noindent b) For $k < L^+$ the local permutation 
is defined by
\bel{perm}
(\pi^{kk+1} (\bfeta))_l = \left\{ \ba{ll} 
\eta_{k+1} & l=k \\
\eta_{k} & l=k+1 \\
\eta_l & \mbox{  else  } \ea \right. .
\ee
For a given $\bfeta$ we use the notation 
$\bfeta^{kk+1} := \pi^{kk+1}(\bfeta)$.
\end{df}

These mappings are invertible and one has
$(\gamma_k)^{-1} = (\gamma_k)^{n}$. 
As an array of occupation variables we can write
\bea 
\gamma_k (\bfeta) & = & (\dots ,\eta_{k-1}, (\eta_k + 1) \mbox{ mod } (n+1), 
\eta_{k+1},  \dots ) \\
\pi^{kk+1}(\bfeta) & = & ( \dots ,\eta_{k-1}, \eta_{k+1}, \eta_k, \eta_{k+2}, 
\dots ) .
\eea

Next we define various functions characterizing the occupation variables
of a configuration $\bfeta$. As a first step we consider a ``lattice''
with only a single site.

\begin{df}
The indicator functions $n^\alpha: \S_{0,n}
\mapsto \{0,1\}$ and $m^\alpha: \S_{0,n}
\mapsto \{0,1\}$ are defined by
\bel{occind1}
n^\alpha(\eta)  :=  \delta_{\eta,\alpha}, \quad m^\alpha(\eta) 
:= \sum_{\beta=\alpha}^n n^\alpha(\eta) %=\Theta(\eta-\alpha+1)
\ee
for $0\leq \alpha \leq n$.
\end{df}

From these we construct the following indicator functions for general 
lattices $\Lambda$.

\begin{df}
The local occupation numbers at site $k\in \Lambda$ of
a configuration $\bfeta \in \S_{0,n}^{\Lambda}$
are defined by
\be 
\label{Def:occind}
n^\alpha_k(\bfeta) :=  n^\alpha(\eta_k), \quad
m^\alpha_k(\bfeta) := \sum_{\beta=\alpha}^n n^\beta_k(\bfeta) 
% = \Theta(\eta_k-\alpha+1)
\ee 
for $0\leq \alpha \leq n$.
\end{df}

By construction $m^{n}_k(\bfeta) = n^{n}_k(\bfeta)$ and $m^{0}_k(\bfeta) =1$.
Notice also that from \eref{Def:occind} one has 
$n^\alpha_k(\bfeta) = \delta_{\eta_k,\alpha}$ which yields the expressions
\bel{partnum}
N^\alpha(\bfeta) = \sum_{k\in\Lambda} n^\alpha_k(\bfeta), \quad
M^\alpha(\bfeta) = \sum_{k\in\Lambda} m^\alpha_k(\bfeta), \quad \alpha\in \S_{0,n}.
\ee
for the particle numbers \eref{Def:partnum}.
We note the trivial, but frequently used identities 
$M^0(\bfeta) = L$, $M^1(\bfeta) = N(\bfeta)$, $M^n(\bfeta)=N^n(\bfeta)$,
and
\be 
\label{occupos}
n^\alpha_k(\bfeta) = \delta_{\eta_k,\alpha} = \sum_{i=1}^{N^\alpha(\bfeta)} 
\delta_{x^\alpha_i(\bfeta),k} .
\ee
Using the representation \eref{sumTheta} in Appendix A 
of the integer Theta-function this identity allows us to write 
\bea
\label{Theta}
\Theta(\eta_k-\eta_l) 
& = & \sum_{\alpha=1}^n \sum_{\beta=0}^{\alpha-1} 
n^\alpha_{k}(\bfeta) n^\beta_{l}(\bfeta) \\
& = & n^0_{l}(\bfeta) (1-n^0_{k}(\bfeta)) +  \sum_{\alpha=2}^n \sum_{\beta=1}^{\alpha-1} 
n^\alpha_{k}(\bfeta) n^\beta_{l}(\bfeta)\\
& = & \sum_{i=1}^{N(\bfx)} \delta_{x_i,k}
\left(1-\sum_{i=1}^{N(\bfx)} \delta_{x_i,l} \right) 
 + \sum_{i=1}^{N(\bfx)} \sum_{j=1}^{N(\bfx)} \delta_{x_i,k} \delta_{x_j,l} 
\, \Theta(\alpha_i-\alpha_j).
\eea
Thus we can express  the two-particle sign function
\bea
\label{sigma}
\sgn(\eta_k-\eta_l) & = & \sum_{\alpha=1}^n \sum_{\beta=0}^{\alpha-1} 
\left(n^\alpha_{k}(\bfeta) n^\beta_{l}(\bfeta) 
- n^\alpha_{l}(\bfeta) n^\beta_{k}(\bfeta)\right) \\
& = & n^0_{l}(\bfeta) - n^0_{k}(\bfeta) + 
\sum_{\alpha=2}^n \sum_{\beta=1}^{\alpha-1} 
\left(n^\alpha_{k}(\bfeta) n^\beta_{l}(\bfeta) 
- n^\alpha_{l}(\bfeta) n^\beta_{k}(\bfeta)\right) \\
& = & \sum_{i=1}^{N(\bfx)} (\delta_{x_i,k} - \delta_{x_i,l} )
+ \sum_{i=1}^{N(\bfx)} 
\sum_{j=1}^{N(\bfx)} \delta_{x_i,k} \delta_{x_j,l} 
\, \sigma(\alpha_i-\alpha_j)
\eea
in terms of occupation numbers $n_k^\alpha$ and also in the
coordinate representation. 

In order to characterize a configuration globally
we introduce the following quantities.

\begin{df} Let $\Lambda$ be one of the four lattices defined 
in \eref{Def:lattice} and $\bfeta\in \S_{0,n}^{\Lambda}$.
The particle balances at site $k$ are defined by
\be   
\label{Def:Nalphafun}
N^{\alpha}_k(\bfeta) := \sum_{l=L^-}^{k-1} n_l^\alpha(\bfeta) 
- \sum_{l=k+1}^{L^+} n_l^\alpha(\bfeta), \quad 
M^{\alpha}_k(\bfeta) := \sum_{\beta=\alpha}^{n} N^{\beta}_k(\bfeta). 
\ee
for $0 \leq \alpha \leq n$.
\end{df}
By construction $M^{n}_k(\bfeta) = N^{n}_k(\bfeta)$.
For $\alpha = 0$ we have 
$M^{0}_k(\bfeta) = 2k-L^+-L^-$.

\subsection{Generator for the $n$-priority ASEP}
\label{Sec:Defprocess}

\subsubsection{The Markov generator and its general matrix formulation}

We recall the definition of a Markov process $\omega_t$ with state space 
$\Omega$ and transition rates $g_{\omega'\omega}$ from a configuration 
$\omega$ to a configuration $\omega'$ 
in terms of a generator $\mathcal{G}$ acting on suitably chosen functions
$f(\omega)$ through the relation
\bel{Markovgen}
\mathcal{G} f(\omega)
= \sum_{\omega' \in \Omega\setminus\omega} 
g_{\omega'\omega} [f(\omega') - f(\omega)].
\ee
The $n$-priority ASEP described informally in the introduction can thus 
be defined as follows.

\begin{df} ($n$-ASEP) 
\label{Def:n-ASEP}
Let $q^{-1} \in (0,1]$ and $\bfeta \in \S_{0,n}^\Lambda$ for 
$\Lambda$ as defined in \eref{Def:lattice} and let 
\be 
\label{localrates}
w^{kk+1}_q(\bfeta) = 
w \left( q^{\sgn(\eta_k-\eta_{k+1})} - \delta_{\eta_k,\eta_{k+1}} \right).
\ee
be the bond hopping rates for the transition rates
\bel{transitionrate}
w_{\bfeta' \bfeta}(q) = \sum_{k=L^-}^{L^+-1} 
w^{kk+1}_q(\bfeta) \delta_{\bfeta',\bfeta^{kk+1}}
\ee
from $\bfeta$ to a locally permutated configuration
$\bfeta^{kk+1}=\pi^{kk+1}(\bfeta)$ defined in \eref{perm}. Then
the $n$-priority ASEP  on $\Lambda$ is the Markov process defined by the generator
\bel{generator}
\mathcal{L} f(\bfeta)
= \sum_{k=L^-}^{L^+-1} \mathcal{L}_{k,k+1} f(\bfeta) 
\ee
with local generators 
\bel{generatorlocal}
\mathcal{L}_{k,k+1} f(\bfeta) = w^{kk+1}_q(\bfeta) [f(\bfeta^{kk+1}) - f(\bfeta)]
\ee
with the convention that 
$L^+ = \infty$ for $\Lambda=\Lambda_+$, $L^- = -\infty$ for $\Lambda=\Lambda_-$,  
and $L^\pm = \pm \infty$ for $\Lambda=\Z$.
\end{df}
The configuration at time $t$ is represented by $\bfeta_t$. 

In the semi-infinite and infinite cases some care needs to be taken
regarding the class of functions to which $f$ belongs as one needs to ensure that
\eref{generator} converges. Following Liggett 
\cite{Ligg85,Ligg99} we note that this is the case for cylinder functions on
$\S_{0,n}^\Z$ (and hence on $\S_{0,n}^\pm)$, 
i.e., functions that depend on only finitely many coordinates.
Going beyond cylinder functions we note the following lemma.

\begin{lmm}
\label{Lemma:ASEPZ}
Let $g$ be a cylinder function on $\S_{0,n}^{\Z}$, and
let $V_r^{>}(\bfeta)):=\sum_{l=r}^\infty (1-n^n_l(\bfeta))$, 
$V_r^{<}(\bfeta)):=\sum_{l=-\infty}^r (1-n^n_l(\bfeta))$
for $r\in\Z$. 
Then for the process defined on $\Lambda=\Z$ one has for $\alpha \geq 1$
and $a\in \R$:\\

\noindent a) $\mathcal{L} f(\bfeta) < \infty$ for
$f(\bfeta)=\rme^{a V_r^{>}(\bfeta)} g(\bfeta)$ and 
$\bfeta \in \mathbb{S}^{\Z}_{n^{>}}$.\\

\noindent b) $\mathcal{L} f(\bfeta) < \infty$ for
$f(\bfeta)=\rme^{a V_r^{<}(\bfeta)} g(\bfeta)$ and 
$\bfeta \in \mathbb{S}^{\Lambda}_{n^{<}}$.
\end{lmm}

\Proof Consider only case (a), case (b) is similar: By definition 
\eref{lfull} we have that $\rme^{a V_r^{>}(\bfeta)} < \infty$ for
any $\bfeta \in \mathbb{S}^{\Z}_{n^{>}}$ which means that $f(\bfeta)$ 
is well-defined on $\mathbb{S}^{\Z}_{n^{>}}$.
We write $f(\bfeta)=\rme^{a V_x^{>}(\bfeta)} 
\rme^{a \sum_{l=k}^{x-1} (1-n^n_l)} g(\bfeta) =: \rme^{a V_x^{>}(\bfeta)} 
\tilde{g}_{x}$ where $x$ is the largest
coordinate on which $g$ depends. 
Since $g$ is a cylinder function one has $x<\infty$ which implies that also
$\tilde{g}_{x}$ is a cylinder function. Because of particle number conservation
one has 
%\be 
$\mathcal{L}_{k,k+1} \rme^{a(n^n_k+n^n_{k+1})} = 0$
%\ee
for any $a\in\R$ 
and therefore $\mathcal{L}_{k,k+1} f(\bfeta) = 0$ for $k>x$. Thus the sum
\eref{generator} with $L^\pm=\pm\infty$ contains only finitely many 
terms and is therefore finite. \qed

Now we focus on the finite lattice $\Lambda_{L}$ and
write the action of the generator \eref{generator} in the so-called quantum 
Hamiltonian form \cite{Sudb95,Schu01}, i.e., in terms of the continuous-time
transition matrix $H$ defined by the matrix elements
\bel{transmatrix}
H_{\bfeta'\bfeta} = \left\{ \ba{ll} 
- w_{\bfeta' \bfeta}(q) \quad & \mbox{ for } \bfeta' \neq \bfeta \\[2mm]
\displaystyle \sum_{\bfeta' \in \S_{0,n}^L\setminus\bfeta} w_{\bfeta' \bfeta}(q) & 
\mbox{ for } \bfeta' = \bfeta .
\ea \right.
\ee
This is a square matrix of dimension 
$d_{n,L}:=|\S_{0,n}^L|=(n+1)^{L}$ where by definition all off-diagonal
elements (the negative transition rates) are non-positive, the diagonal
elements are all non-negative and
in each column the sum of all matrix elements is equal to 0, expressing
probability conservation $\mathcal{L} f =0$ for the identity function
$f(\bfeta)=1$.

In terms of the matrix elements $H_{\bfeta'\bfeta}$
the defining equation \eref{Markovgen} then has the form
\bel{generatorH}
\mathcal{L} f(\bfeta) = - \sum_{\bfeta' \in \S_{0,n}^L} f(\bfeta') H_{\bfeta'\bfeta}
\ee
which can be interpreted as a matrix multiplication $-\mathbf{f} H$ where
according to standard convention for matrix multiplication $\mathbf{f}$ is 
understood as a row vector with entries $f(\bfeta)$.

With \eref{generator}, \eref{generatorlocal} we can write \eref{generatorH}
as
\bel{generatorH2}
\mathcal{L} f(\bfeta) = - \sum_{k=L^-}^{L^+-1} \sum_{\bfeta' \in \S_{0,n}^L} 
f(\bfeta') (h_{k,k+1})_{\bfeta'\bfeta}
\ee
with the local hopping matrices $h_{k,k+1}$ which are the continuous-time
transition matrices of the process restricted to the bond $(k,k+1)$.
In slight abuse of language we shall call also 
\bel{H}
H =  \sum_{k=L^-}^{L^+-1} h_{k,k+1}
\ee
the generator of the process. 

\subsubsection{The tensor basis}

In order to write $H$ explicitly it is natural to choose the canonical
basis $\mathfrak{B}=\{\mathfrak{b}(i),1\leq i \leq d_{n,L}\}$ which spans the 
vector space $\C^{d_{n,L}}$.
The sum \eref{generatorH} does not uniquely define the 
matrix $H$ as the arrangement of the elements $f(\bfeta)$ in a vector
$\mathbf{f}$ is not specified by this sum. One has still the freedom
to define the mapping $\iota: \S_{0,n}^L \mapsto \{1,\dots,d_{n,L}\}$
that specifies which canonical basis vector 
$\mathfrak{b}(\iota(\bfeta))$ corresponds to a given configuration $\bfeta$.
We use the natural quantum Hamiltonian form \cite{Sudb95,Schu01}
where the ordering of the basis is given by the numerical $n$-ary representation of a 
configuration $\bfeta$ defined
as follows.
\begin{df} (Basis order)
\label{Def:basis}
Let $\C^{d_{n,L}}$ be the 
$d_{n,L}$-dimensional vector space over $\C$ with canonical basis vectors
$\mathfrak{b}(i) = (0,\dots,0,1,0,\dots,0)$ in row form with entry 1 for
component $i$ with $1\leq i \leq d_{n,L}$ and zero else. 
The basis $\mathfrak{E}=\{\bra{\bfeta}, \bfeta\in\S_{0,n}^L\}$ 
of $\C^{d_{n,L}}$ is defined by the row vectors $\bra{\bfeta}
= \mathfrak{b}(\iota(\bfeta))$ with
\bel{ternary}
\iota(\bfeta) = 1 + \sum_{k=L^-}^{L^+} \eta_k \, (n+1)^{k-L^-}.
\ee
The dual basis of column vectors is given by $\ket{\bfeta}:=\bra{\bfeta}^T$
where the superscript $T$ denotes transposition.
\end{df}

The term ``quantum Hamiltonian formalism'', which also motivates the use
of the bra ($\bra{\cdot}$) and ket ($\ket{\cdot}$) symbols for vectors, 
will become clear below.
A general row vector with entries $f(\bfeta)$ is denoted by $\bra{f}$ and a
column vector with entries $g(\bfeta)$ is denoted by $\ket{g}$. 
The ordering of the basis induced by the numerical presentation of the 
configurations induces a tensor structure which is defined by the Kronecker 
product \eref{Def:Kronecker} defined in the Appendix.
Interpreting a row vector of dimension $d$ as a  $1 \times d$-matrix we 
then have:

\begin{prop} (Tensor basis)
For $\eta \in  \S_{0,n}$ let $(\eta|$ be the canonical
$(n+1)$-dimensional basis vectors of
$\C^{n+1}$ with component $1$ at position $1+\eta$ 
and zero else and let $|\eta)=(\eta|^T$ be the dual basis vector. 
Then one has
\be 
\label{tensorbasis}
\bra{\bfeta} = (\eta_{L^-}| \otimes (\eta_{L^-+1}| \otimes \dots \otimes 
(\eta_{L^+}| , \quad
\ket{\bfeta} =
|\eta_{L^-}) \otimes |\eta_{L^-+1}) \otimes \dots \otimes |\eta_{L^+})
\ee
with $\bra{\bfeta}$ and $\ket{\bfeta}$ as in Definition \eref{Def:basis}.
\end{prop}
\Proof This is an immediate consequence of the definition \eref{Def:Kronecker}
of the Kronecker product and the $n$-ary representation \eref{ternary} of a 
configuration $\bfeta$.
\qed

Following quantum mechanical convention we omit the tensor symbol $\otimes$
in the Kronecker product of bra and ket vectors. In particular, we use
\bel{tensprod}
\ket{g}\bra{f} \equiv \ket{g} \otimes \bra{f} .
\ee
By Definition \eref{Def:Kronecker} this is $d_{n,L} \times d_{n,L}$-matrix 
$C$ with matrix elements 
$C_{i,j}=g_i f_j$. Specifically, we have the representation
\bel{unitmatrix}
\boldsymbol{1} = \sum_{\bfeta \in \S_{0,n}^L} \ket{\bfeta}\bra{\bfeta}.
\ee
of the $d_{n,L}$-dimensional unit matrix and
\bel{projectormatrix}
\hat{\wp}^\alpha = \sum_{\beta=\alpha}^n \sum_{\bfeta \in \S_{N^\alpha}^L}
\ket{\bfeta}\bra{\bfeta}
\ee
of the projector matrix with the property 
$\bra{\bfeta} \hat{\wp}^\alpha = \bra{\bfeta} \wp^\alpha_{\vec{N}}(\bfeta)$
derived from the projector definition \eref{particleproj}.

In addition to the tensor product we need an inner product from $\C^{d_{n,L}}$
and its dual to $\C$.

\begin{df}
The inner product for row vectors $\bra{f} = \sum_{\bfeta} f(\bfeta) \bra{\bfeta}$
and column vectors $\ket{g}= \sum_{\bfeta} g(\bfeta) \ket{\bfeta}$ is defined by
\bel{inprod}
\inprod{f}{g} := \sum_{\bfeta \in \S_{0,n}^L} f(\bfeta) g(\bfeta) .
%= \sum_{i=1}^{d_{n,L}} f(i)g(i).
\ee
\end{df}
This implies the biorthogonality relation $\inprod{\bfeta'}{\bfeta} = 
\delta_{\bfeta',\bfeta}$. Next we state some general rules of multilinear algebra 
that motivate the omission of the tensor symbol in \eref{tensprod},
highlight a factorization property of the tensor basis under the inner product, and
illustrates the use of the representation \eref{unitmatrix} of the unit matrix.

\begin{lmm}
\label{Lem:inprodprops}
a) Let $C=\ket{g}\bra{f}$ be a tensor matrix according to definition \eref{tensprod}.
Then the inner product with vectors $\bra{a}$ and $\ket{b}$ is given by
\be 
\bra{a} C \ket{b} \equiv \bra{a} (\ket{g}\otimes \bra{f}) \ket{b} 
\equiv \bra{a} (\ket{g}\bra{f}) \ket{b} 
= \inprod{a}{g} \cdot \inprod{f}{b}
\ee
where $\cdot$ denotes ordinary multiplication in $\C$.\\

\noindent b) Let the vectors 
$\bra{f} = (f_L^-|\otimes(f_{L^-+1}|\otimes \dots \otimes (f_{L^+}|$
and $\ket{g} = (g_L^-|\otimes(g_{L^-+1}|\otimes \dots \otimes (g_{L^+}|$ be tensor
products. Then
\bel{inprodfac} 
\inprod{f}{g} = \prod_{k=L^-}^{L^+} (f_k|g_k) = 
\prod_{k=L^-}^{L^+} \left( \sum_{\alpha=0}^n f_k(\alpha) g_k(\alpha) \right).
\ee

\noindent c) For any pair of functions $f$ and $g$ the inner product can be 
expanded in a complete basis as 
\be 
\inprod{f}{g} = \sum_{\bfeta\in \S_{0,n}^L} \inprod{f}{\bfeta}\inprod{\bfeta}{g}.
\ee
\end{lmm}

\begin{rem} 
The expectation $\E^\mu f := \sum_{\bfeta} f(\bfeta) \mu(\bfeta)$ 
of a function $f(\bfeta)$ under a 
probability measure $\mu(\bfeta)$ can be written as the inner product
$\E^\mu f = \inprod{f}{\mu}$.
Defining the summation vector 
\bel{sumvec}
\bra{s} := \sum_{\bfeta \in \S_{0,n}^L} \bra{\bfeta}
\ee
where all entries are equal to 1 and for a function $f(\bfeta)$ 
the diagonal matrix
\bel{diagf}
\hat{f} := \sum_{\bfeta \in \S_{0,n}^L} f(\bfeta) \ket{\bfeta}\bra{\bfeta}
\ee
one can write $\bra{f} = \bra{s} \hat{f}$ and therefore
$\E^\mu f = \bra{s} \hat{f} \ket{\mu}$.
\end{rem}
This remark highlights the role
of the representation of a function $f$ as a diagonal matrix which
we shall generally denote by the
circumflex ( $\hat{}$ )-symbol. 

We also use a diagonal matrix representation of a probability measure.
\begin{prop} 
\label{revprocess}
Let $\mu >0$ be a strictly positive
reversible measure for a process with generator $H$. Then with 
the diagonal matrix representation 
\bel{revmeasmatrix}
\hat{\mu} :=  \sum_{\bfeta} \mu(\bfeta)  \ket{\bfeta}  \bra{\bfeta} 
\ee
of $\mu$ the transformation property
\bel{detbalH}
\hat{\mu}^{-1} H \hat{\mu} = H^T
\ee
is the condition of detailed balance.
\end{prop}

\Proof ``Sandwiching'' \eref{detbalH} with $\bra{\bfeta'}$ and $\ket{\bfeta}$
one finds
\bel{detbal}
\mu(\bfeta') w_{\bfeta'\bfeta} =  \mu(\bfeta) w_{\bfeta\bfeta'}.
\ee
which is indeed the detailed-balance condition for reversible measures.
\qed

\subsubsection{Explicit construction of $H$}

In the quantum Hamiltonian formalism
the functions of configurations defined in \eref{Def:functions}
turn into endomorphisms on the vector space $\C^{d_{n,L}}$,
represented by matrices.
In order to construct $H$ explicitly 
we first define for $1\leq \alpha \leq n$
the following single-site matrices of dimension $n+1$:

\begin{df}
\label{Def:singlesite}
For $1\leq \alpha \leq n$ the single-site raising and lowering operators 
$\sigma^{\alpha,\pm}$ and for $0\leq \alpha \leq n$ the single-site projectors 
$\hat{n}^\alpha$ are defined by
\bel{fundgln}
\sigma^{\alpha,+} := |\alpha-1)(\alpha|, \quad \sigma^{\alpha,-} := |\alpha)
(\alpha-1|, \quad \hat{n}^\alpha = |\alpha)(\alpha|.
\ee
The cyclic flip operator $\gamma$ and the species flip $\sigma^{\alpha\beta}$
are defined by
\bel{cyclic}
\gamma := |n)(0| + \sum_{\alpha=1}^n \sigma^{\alpha,+}, \quad 
\sigma^{\alpha\beta} := |\alpha)(\beta|.
\ee
%the inversion operator
%\bel{inversion}
%\chi := \sum_{\alpha=0}^n |\alpha)(n-\alpha|
%\ee
\end{df}
Notice the cyclic property $\gamma^{n} = \gamma^{-1}$ and the 
representation 
\bel{gammanlocal}
\sigma^{\alpha\beta} = \hat{n}^\alpha \gamma^{\beta-\alpha}= 
\gamma^{\beta-\alpha}\hat{n}^\beta
\ee
of the species flip operator.
We denote the unit matrix of dimension $n+1$ by $\mathds{1}$.

From these matrices we construct local operators acting non-trivially on the 
configuration at site $k \in \Lambda_{L}$ as follows.

\begin{df}
\label{Def:localop}
Let $a$ be matrix of dimension $n+1$. For $k\in\Lambda^L$ the local operator 
$a_k$ of dimension $d_{n,L}$ is defined by
\bel{localop}
a_k = \mathds{1}^{\otimes (k-L^-)} \otimes a \otimes \mathds{1}^{\otimes (L^+-k)}
\ee
with the convention that $\mathds{1}^{\otimes 0} = 1$.
\end{df}

Notice that for any pair of matrices $a,b$ one has the 
commutator property 
\bel{localcommutator}
\comm{a_k}{a_l}=0 
\ee 
and for any pair of tensor products $A=\prod_{k\in\Lambda} a_k$ and
$B=\prod_{k\in\Lambda} b_k$ one has the tensor factorization property
\bel{tensorfactor}
AB = \prod_{k\in\Lambda} (ab)_k.
\ee
Both these properties, which arise from the multilinearity of the tensor
product, will be used throughout this work.

With this construction we obtain local operators $\gamma_k$
with the properties $\bra{\bfeta} \gamma_k = \bra{\bfeta_k^{+}}$ and
$\gamma_k \ket{\bfeta} = \ket{\bfeta_k^{-}}$. The global cyclic
flip operator is then given by
\bel{globcycl}
\Gamma := \prod_{k=L^-}^{L^+} \gamma_k
\ee
and has the properties $\bra{\bfeta} \Gamma = \bra{\bfeta^{+}}$ and
$\Gamma \ket{\bfeta} = \ket{\bfeta^{-}}$.

We are now in a position to write the generator $H$ \eref{H}
in explicit form, using the fact that the inner product implies 
$H_{\bfeta'\bfeta} = \bra{\bfeta'} H \ket{\bfeta}$.
%, and state one of its symmetries.

\begin{prop}
Define
the single-bond hopping matrices
\bel{singlebond}
h_{k,k+1} := - w \sum_{\alpha=1}^n \sum_{\beta=0}^{\alpha-1} 
\left( q \sigma^{\beta \alpha}_k \sigma^{\alpha \beta}_{k+1} +
q^{-1} \sigma^{\alpha \beta}_k \sigma^{\beta \alpha}_{k+1}\right) + \hat{w}_{k,k+1}
\ee
where 
\be 
\hat{w}_{k,k+1} = w \sum_{\alpha=1}^n \sum_{\beta=0}^{\alpha-1} 
\left( q \hat{n}^\alpha_k \hat{n}^\beta_{k+1} + 
q^{-1} \hat{n}^\beta_{k} \hat{n}^\alpha_{k+1} \right).
\ee
The generator $H$ of the $n$-priority ASEP on $\Lambda_{L}$ 
defined by \eref{generator} and \eref{generatorlocal} is given in quantum
Hamiltonian form by the matrix \eref{H} with local transition matrices
\eref{singlebond}.
\end{prop}

\Proof Since $\hat{n}_k^\alpha \ket{\bfeta} 
= n_k^\alpha(\bfeta) \ket{\bfeta}$ the diagonal part, which
yields the negative contribution to \eref{generatorlocal}, follows from the
expression \eref{sigma} of the sign-function in terms of the single-site
projectors and the projector Lemma \eref{projlem}.
For the off-diagonal part notice that 
\be 
\bra{\bfeta'} \gamma^{-1}_k \ket{\bfeta} = \delta_{\bfeta',\bfeta^{k,+}}
\ee
with $\bfeta^{k,+} = \gamma_k(\bfeta)$ defined in \eref{flip} and therefore
\be 
\bra{\bfeta'} \sigma^{\alpha \beta}_k \ket{\bfeta} = 
n_k^\alpha(\bfeta')  \delta_{\bfeta',\bfeta^{k,\alpha-\beta}}
= \delta_{\eta_k',\alpha} \delta_{\eta_k,\beta}.
\ee
This yields
\be 
\bra{\bfeta'} \sigma^{\alpha \beta}_k \sigma^{\beta \alpha}_{k+1} 
\ket{\bfeta} =
\delta_{\eta_k',\alpha} \delta_{\eta_k,\beta}
\delta_{\eta_{k+1},\alpha} \delta_{\eta_{k+1}',\beta}
= n_{k}^\beta(\bfeta) n_{k+1}^\alpha(\bfeta) 
 \delta_{\bfeta',\bfeta^{kk+1}} 
\ee
with the local permutation $\bfeta^{kk+1}$ \eref{perm}.
Taking the summation over $\alpha$ and $\beta$ and using \eref{Theta}
yields the offdiagonal part
of \eref{generatorlocal}. \qed

\begin{rem}
The generator $H$ has a structure reminiscent of the quantum
Hamiltonian of the Perk-Schultz quantum chain \cite{Perk81}
\bel{HPS}
H^{PS} = \sum_{k=L^-}^{L^+-1} h^{PS}_{k,k+1}
\ee
with the single-bond matrices
\bel{singlebondPS}
h^{PS}_{k,k+1} = - w \sum_{\alpha=1}^n \sum_{\beta=0}^{\alpha-1} 
\left( \sigma^{\beta,\alpha}_k \sigma^{\alpha,\beta}_{k+1} +
\sigma^{\alpha,\beta}_k \sigma^{\beta,\alpha}_{k+1}\right) + \hat{w}_{k,k+1}.
\ee
Notice that unlike the generator \eref{H} the Hamiltonian $H^{PS}$ is symmetric 
(and hence Hermitian as should be the case for a
quantum system).
\end{rem}

\subsection{The quantum algebra $U_q[\mathfrak{gl}(n+1)]$}

The Perk-Schultz quantum chain is an integrable model which brings in 
the notion of quantum algebras. We first introduce the 
quantum algebra $U_q[\mathfrak{gl}(n+1)]$ \cite{Jimb85,Jimb86}
in terms of abstract generators and then give representation
matrices that satisfy its defining relations.

\begin{df}
\label{Def:Uqglnp1}
The quantum algebra $U_q[\mathfrak{gl}(n+1)]$ \cite{Jimb85,Jimb86} is 
the associative algebra over $\C$
generated by $\mathbf{L}_\alpha^{\pm 1}$, $\alpha=0,\dots,n$ and 
$\mathbf{X}^\pm_\alpha$, $\alpha=1,\dots,n$ and unit $I$ with relations 
\bea
\label{Uqglndef1}
& & \mathbf{L}_\alpha^{\pm 1} \mathbf{L}_\alpha^{\mp 1} = I  \\
\label{Uqglndef2}
& & \comm{\mathbf{L}_\alpha}{\mathbf{L}_\beta} = 0  \\
\label{Uqglncomm2}
& & \mathbf{L}_\alpha \mathbf{X}^\pm_\beta = q^{\pm (\delta_{\alpha,\beta-1} - 
\delta_{\alpha,\beta})/2} \mathbf{X}^\pm_\beta \mathbf{L}_\alpha\\
\label{Uqglncomm3}
& & \comm{\mathbf{X}^+_\alpha}{\mathbf{X}^-_\beta} = 
\delta_{\alpha,\beta} \frac{(\mathbf{L}_{\alpha-1}\mathbf{L}_\alpha^{-1})^2 - 
(\mathbf{L}_{\alpha-1}\mathbf{L}_\alpha^{-1})^{-2}}{q-q^{-1}}\\
\label{UqglnSerre1}
& & \comm{\mathbf{X}^\pm_\alpha}{\mathbf{X}^\pm_\beta} 
= 0 \quad |\alpha-\beta| \neq 1, \\
\label{UqglnSerre2}
& & (\mathbf{X}^\pm_\alpha)^2 \mathbf{X}^\pm_\beta 
- [2]_q \mathbf{X}^\pm_\alpha \mathbf{X}^\pm_\beta \mathbf{X}^\pm_\alpha + 
\mathbf{X}^\pm_\beta (\mathbf{X}^\pm_\alpha)^2 = 0 \quad |\alpha-\beta| = 1.
\eea
\end{df}

\subsubsection{Fundamental representation of $\mathfrak{gl}(n+1)$}

We define for $0 \leq \alpha \leq n$
\bel{Def:Nalpha}
\mathbf{L}_\alpha = q^{\frac{1}{2}\mathbf{N}_\alpha}.
\ee
In terms of the $\mathbf{N}_\alpha$ the defining relations of the Lie algebra 
$\mathfrak{gl}(n+1)$
are given by the limit $q\to 1$ of \eref{Uqglndef2} - \eref{UqglnSerre2}:
\bea
\label{glndef1}
& & \comm{\mathbf{N}_\alpha}{\mathbf{N}_\beta} = 0  \\
\label{glncomm2}
& & 
\comm{\mathbf{N}_\alpha}{\mathbf{X}^\pm_\beta} = 
\pm (\delta_{\alpha,\beta-1} - \delta_{\alpha,\beta} ) \mathbf{X}^\pm_\beta \\
\label{glncomm3}
& & \comm{\mathbf{X}^+_\alpha}{\mathbf{X}^-_\beta} = \delta_{\alpha,\beta}
\left( \mathbf{N}_{\alpha-1}
- \mathbf{N}_{\alpha} \right) \\
\label{glnSerre1}
& & \comm{\mathbf{X}^\pm_\alpha}{\mathbf{X}^\pm_\beta} 
= 0 \quad |\alpha-\beta| \neq 1, \\
\label{glnSerre2}
& & \comm{\mathbf{X}^\pm_\alpha}{\comm{\mathbf{X}^\pm_\alpha}{\mathbf{X}^\pm_\beta}} 
= 0 \quad |\alpha-\beta| = 1.
\eea

It is well-known (and easy to verify) that the matrices \eref{fundgln} form the 
fundamental representation of the Lie algebra $\mathfrak{gl}(n+1)$ 
\eref{glndef1} - \eref{glnSerre2} via the algebra homomorphism 
$\mathbf{X}^\pm_\alpha \mapsto 
\sigma^{\alpha,\pm}$, $\mathbf{N}_\alpha \mapsto \hat{n}^{\alpha}$.
It has been pointed out \cite{Jimb86} that with \eref{Def:Nalpha}
these matrices then also form a representation of the quantum algebra 
$U_q[\mathfrak{gl}(n+1)]$ \eref{Uqglndef1} - \eref{UqglnSerre2}.
We remark that $\sigma^{\alpha,\pm}$ are nilpotent of degree 2, i.e.,
$(\sigma^{\alpha,\pm})^2=0$.

\subsubsection{Relation between $U_q[\mathfrak{gl}(n+1)]$ and 
$U_q[\mathfrak{sl}(n+1)]$}

\begin{df}
\label{Def:Uqslnp1}
Let $A$ be the Cartan matrix of simple Lie algebras
of type $A_{n+1}$ with matrix elements
$A_{\alpha\beta} = 2 \delta_{\alpha,\beta} - 
\delta_{\alpha,\beta-1} \delta_{\alpha,\beta+1}$
and 
\bel{complement}
\mathbf{H}_\alpha = \mathbf{N}_{\alpha-1} - \mathbf{N}_{\alpha}, \quad 
1 \leq \alpha \leq n .
\ee
Then the quantum algebra $U_q[\mathfrak{sl}(n+1)]$ is the subalgebra 
of $U_q[\mathfrak{gl}(n+1)]$ generated by 
$q^{\pm \mathbf{H}_\alpha/2}$ and $\mathbf{X}^\pm_\alpha$ with relations
\bea
& & q^{\mathbf{H}_\alpha/2} q^{-\mathbf{H}_\alpha/2} 
= q^{-\mathbf{H}_\alpha/2} q^{\mathbf{H}_\alpha/2} = I\\
\label{Uqslncomm1b}
& & q^{\mathbf{H}_\alpha/2} q^{\mathbf{H}_\beta/2} 
= q^{\mathbf{H}_\beta/2} q^{\mathbf{H}_\alpha/2} \\
\label{Uqslncomm2b}
& & q^{\mathbf{H}_\alpha} \mathbf{X}^\pm_\beta q^{-\mathbf{H}_\alpha} 
= q^{\pm A_{\alpha \beta}} \mathbf{X}^\pm_\beta\\
\label{Uqslncomm3b}
& & \comm{\mathbf{X}^+_\alpha}{\mathbf{X}^-_\beta} 
= \delta_{\alpha\beta} [\mathbf{H}_\alpha]_q.
\eea
and \eref{UqglnSerre1}, \eref{UqglnSerre2}.
\end{df}
The fact that $U_q[\mathfrak{sl}(n+1)]$ is a subalgebra of 
$U_q[\mathfrak{gl}(n+1)]$ can be seen by noticing that 
$\prod_{\alpha=0}^n \mathbf{L}_\alpha$ belongs to the center of 
$U_q[\mathfrak{gl}(n+1)]$. The fundamental representation
of both $\mathfrak{sl}(n+1)$ and $U_q[\mathfrak{sl}(n+1)]$ is formed by the
set of matrices $\sigma^{\alpha,\pm}$ defined in \eref{fundgln} and
\bel{Def:halpha}
\hat{h}^\alpha = \hat{n}^{\alpha-1} - \hat{n}^\alpha
\ee
with $1\leq \alpha \leq n$.

\subsubsection{Coproduct representation of $U_q[\mathfrak{gl}(n+1)]$}

The coproduct is an algebra homomorphism defined by \cite{Jimb85}
\bel{coprod}
\Delta(\sigma^{\alpha,\pm}) = \sigma^{\alpha,\pm} \otimes 
q^{-\frac{1}{2}\hat{h}^\alpha} + q^{\frac{1}{2}\hat{h}^\alpha} \otimes 
\sigma^{\alpha,\pm}, \quad
\Delta(\hat{n}^\alpha) = \hat{n}^\alpha \otimes \mathds{1} + \mathds{1} 
\otimes \hat{n}^\alpha.
\ee
Iteration yields the representation matrices
\bea
\label{coprodrepX}
X^{\pm}_\alpha & = & \sum_{k=L^-}^{L^+} X^{\alpha,\pm}(k) \\
\label{coprodrepN}
\hat{N}^\alpha & = & \sum_{k=L^-}^{L^+} \hat{n}^\alpha_k
\eea
with
\bea
X^{\alpha,\pm}(k) & = & \left(q^{\frac{1}{2} 
\hat{h}^\alpha} \right)^{\otimes (k-L^-)} \otimes \sigma^{\alpha,\pm} \otimes
\left(q^{-\frac{1}{2} \hat{h}^\alpha} \right)^{\otimes (L^+-k)} \\
& = & q^{\frac{1}{2} \sum_{l=L^-}^{k-1} \hat{h}^\alpha_l} 
\sigma^{\alpha,\pm}_k q^{- \frac{1}{2} \sum_{l=k+1}^{L^+} \hat{h}^\alpha_l} 
\eea
where
\bea 
\hat{n}^\alpha_k  & = & \mathds{1}^{\otimes (k-L^-)} \otimes \hat{n}^\alpha \otimes
\mathds{1}^{\otimes (L^+-k)} \\
\sigma^{\alpha,\pm}_k & = & \mathds{1}^{\otimes (k-L^-)} \otimes 
\sigma^{\alpha,\pm} \otimes \mathds{1}^{\otimes (L^+-k)}.
\eea
The unit $I$ is represented by the $d_{n,L}$-dimensional unit matrix
\bel{rep1}
\mathbf{1} := \mathds{1}^{\otimes L}.
\ee

The crucial property of this representation are the commutation relations
\bel{Hsym}
\comm{H^{PS}}{X^{\alpha,\pm}} = \comm{H^{PS}}{N^\alpha} = 0
\ee
which express the symmetry of the Perk-Schultz Hamiltonian $H^{PS}$ \eref{HPS}
under the action of the quantum algebra $U_q[\mathfrak{gl}(n+1)]$ \cite{Alca93}.

Notice that for the symmetric case $q=1$ we shall denote the representations
$X^{\pm}_\alpha$ by $S^{\pm}_\alpha$. Together with the diagonal matrices
$N^\alpha$ They form a tensor representation of
$\mathfrak{gl}(n+1)$ as defined in \eref{glndef1} - \eref{glnSerre2}.

\subsection{Remarks on the hydrodynamic limit}

We prove some simple results on locally conserved currents 
that suggest certain properties of the hydrodynamic limit of the 
$n$-priority ASEP.

\begin{prop}
Let $\mathcal{L}$ be the generator of the $n$-species priority ASEP on the finite
lattice $\Lambda_{L}$. The local indicators $n^\alpha_{k}(\bfeta)$ 
\eref{Def:occind} satisfy the discrete 
continuity equation
\bel{continuity}
\mathcal{L} n^\alpha_k = 
j^\alpha_{k-1} - j^\alpha_k, \quad L^- \leq k \leq L^+
\ee
with the locally conserved instantaneous currents
\bel{conscurr}
j^\alpha_k = w \left( \sum_{\beta=0}^{\alpha-1} 
\left(q n^\alpha_{k} n^\beta_{k+1} - q^{-1} n^\beta_k n^\alpha_{k+1}\right)
- \sum_{\beta=\alpha+1}^{n} \left(q n^\beta_k n^\alpha_{k+1} 
- q^{-1} n^\alpha_{k} n^\beta_{k+1}\right)
\right).
\ee
for $L^- \leq k \leq L^+$ and $j^\alpha_{L^--1}=j^\alpha_{L^++1}=0$ for all
$\alpha \in \S_{0,n}$.
\end{prop}

\Proof
In order to compute $\mathcal{L} n_k^\alpha(\bfeta)$ we note that because of the
projector property
\bea
\label{aux1a}
w^{kk+1}_q(\bfeta) n^\alpha_{k}(\bfeta) & = & w n^\alpha_{k}(\bfeta) 
\left( q \sum_{\beta=0}^{\alpha-1} n^\beta_{k+1}(\bfeta)
+ q^{-1} \sum_{\beta=\alpha+1}^{n} n^\beta_{k+1}(\bfeta) \right) \\
\label{aux1b}
w^{kk+1}_q(\bfeta) n^\alpha_{k+1}(\bfeta) & = & w n^\alpha_{k+1}(\bfeta) 
\left( q^{-1} \sum_{\beta=0}^{\alpha-1} n^\beta_k(\bfeta)
+ q \sum_{\beta=\alpha+1}^{n} n^\beta_k(\bfeta)\right).
\eea
Moreover, from \eref{perm} one gets
\bea
n_k^\alpha(\bfeta^{rr+1}) & = & n_k^\alpha(\bfeta)  
+ (n_{k-1}^\alpha(\bfeta) - n_{k}^\alpha(\bfeta)) \delta_{k,r+1}
- (n_{k}^\alpha(\bfeta) - n_{k+1}^\alpha(\bfeta)) \delta_{k,r} \nonumber \\
& = & n_k^\alpha(\bfeta) + 
(n_{r}^\alpha(\bfeta) - n_{r+1}^\alpha(\bfeta)) (\delta_{k,r+1}-\delta_{k,r}).
\eea
which yields the locally conserved instantaneous currents 
$j^\alpha_k(\bfeta) = - w^{kk+1}_q(\bfeta)(n^\alpha_{k+1}(\bfeta) - n^\alpha_k(\bfeta))$
for which we obtain from \eref{aux1a}, \eref{aux1b} 
the explicit expressions \eref{conscurr}. \qed

We can alternatively write the instantaneous currents in the form
\bea
j^\alpha_k 
& = &  w \left( q n^\alpha_{k}  - q^{-1}  n^\alpha_{k+1}\right) 
- w \left( q n^\alpha_{k} n^\alpha_{k+1} 
- q^{-1} n^\alpha_k n^\alpha_{k+1}\right) \nonumber \\
& & - w \left(q - q^{-1}\right) 
\sum_{\beta=\alpha+1}^{n} \left(n^\alpha_{k} n^\beta_{k+1} 
- n^\beta_k n^\alpha_{k+1}\right).
\eea
For $q=1$ this reduces to the linear diffusive current 
$j^\alpha_k = n^\alpha_{k}  - n^\alpha_{k+1}$ 
which means that the symmetric $n$-component exclusion process is a 
decoupled gradient system. For $q\neq 1$ a naive
guess (based on the coupling to indicators $\beta\neq\alpha$ in the
expressions for the microscopic currents $j^\alpha$) seems to suggest
that the hydrodynamic limit of the model would correspond to coupled
Burgers equations similar to those treated in \cite{Spoh14}. However, 
as demonstrated by the next result, these equations can be decoupled 
very simply.

\begin{prop}
Let $\mathcal{L}$ be the generator of the $n$-species priority ASEP on the finite
lattice $\Lambda_{L}$. The local indicators $m^\alpha_{k}(\bfeta)$ 
\eref{Def:occind} satisfy the discrete 
continuity equation
\bel{continuity2}
\mathcal{L} m^\alpha_k = 
\tilde{j}^\alpha_{k-1} - \tilde{j}^\alpha_k , \quad L^- \leq k \leq L^+
\ee
with the locally conserved instantaneous currents
\bel{conscurr2}
\tilde{j}^\alpha_k = w \left(q m^\alpha_k(1-m^\alpha_{k+1}) - q^{-1} (1-m^\alpha_k)m^\alpha_{k+1} \right)
\ee
for $L^- \leq k \leq L^+$ and $j^\alpha_{L^--1}=j^\alpha_{L^++1}=0$ for all
$\alpha \in \S_{0,n}$.
\end{prop}

\Proof
Straightforward computation yields 
$\tilde{j}^\alpha_k = - w^{kk+1}_q(m^\alpha_{k+1} - m^\alpha_k)$
and after some algebra involving reshuffling of indices one gets from 
\eref{conscurr} the expression \eref{conscurr2}. \qed

\begin{rem}
\label{redproc}
Even for $q\neq 1$ there is no cross-coupling between different indices 
$\alpha,\beta$ in the currents \eref{conscurr2}. This stems from the fact 
that the process for 
the occupation variables $m^\alpha_k$ is the same as for a 
first-class particle $n$ that sees all other particles as
vacancies and thus does not depend on $\alpha$. 
\end{rem}

This fact provides some intuition on the hydrodynamic
limit \cite{Spoh91,Kipn99} of the process. We shall assume that
each species $\alpha$ is represented by $\rho^\alpha:=N^\alpha/L$ particles such 
that in the hydrodynamic limit $L\to\infty$ the densities $\rho^\alpha$ remain 
non-zero. For our purposes an informal discussion is sufficient.

The Markov projection on the microscopic occupation numbers $m_k$ and the
resulting decoupling \eref{continuity2} shows that on macroscopic scale 
the $n$-species the priority ASEP with $q\neq 1$ is governed by a system of 
{\it decoupled} inviscid Burgers equations 
\bel{Burgers} 
\partial_t \sigma^\alpha(x,t) + w(q-q^{-1})\partial_x \left[\sigma^\alpha(x,t)
(1-\sigma^\alpha(x,t))\right] = 0
\ee
for the local densities 
$\sigma^\alpha(x,t) := \sum_{\beta=\alpha}^n \rho^\alpha(x,t)$,
in complete analogy to the usual ASEP \cite{Reza91}. 
In infinite volume the Riemann problem can be solved with
the method of characteristics and shock stability \cite{Lax73,Lax06}. 
In particular, for any fixed $\alpha$ there exists a {\it shock} solution
with density
\bel{hydrodensityprofile}
\sigma^\alpha(x,t) = \left\{ \ba{ll}
\sigma^\alpha_- & x < x^\alpha_s(t) \\ 
\sigma^\alpha_+ &  x > x^\alpha_s(t)
\ea \right.
\ee
and deterministically moving shock position 
\bel{shockpos}
x_s^\alpha(t) = x_s^\alpha(0) + v^\alpha_s t
\ee 
with shock velocity
\bel{shockvel}
v^\alpha_s = w(q-q^{-1}) (1-\sigma^\alpha_+ - \sigma^\alpha_-)
\ee
arising from the Rankine-Hugoniot condition \cite{Lax73}
\bel{RH}
v_s(c_+,c_-) = \frac{j_+ - j_-}{c_+-c_-}
\ee
for discontinuities from $c_-$ to $c_+$ and currents $j_\pm$
in the two branches of the shock, which in the present case are given by
$j^\alpha_{\pm} = w(q-q^{-1}) \sigma^\alpha_{\pm}(1-\sigma^\alpha_{\pm})$
for the shock densities $\sigma^\alpha_{\pm}$.
Since by \eref{shockvel} 
$v^\alpha_s-v^\beta_s>0$ for $\alpha>\beta$
the shock positions for the individual modes $\alpha$ satisfy
$x_s^n(t) > x_s^{n-1}(t) > \dots > x_s^1(t)$ at all times $t\geq 0$.

More complex weak solutions allow for consecutive shocks in each mode
at positions $x_s^{\alpha,i}(t)$. According to \eref{shockvel}
the corresponding shock velocities satisfy $v^{\alpha,i}_s > v^{\alpha,i+1}_s$
for all $i$, so that after a finite macroscopic time neighboring shocks
coalesce \cite{Ferr00} and eventually only a single shock of type 
$\alpha$ remains. 
The question then arises what these macroscopic
discontinuities look like on microscopic scale, see \cite{Derr93} for the stationary
limit of a single travelling shock in the standard ASEP and 
for the time-dependent case see 
\cite{Beli02,Bala10}. This issue will be further addressed below.

For the finite system with reflecting boundaries we 
scale the lattice edges $L^\pm = x^\pm L$ such that for $L\to\infty$
macroscopic volume $\ell = x^+ - x^-$ remains finite.
The reflecting boundary
conditions correspond on macroscopic scale to the evolution of 
$\sigma^\alpha(x,t)$ on the interval $[x^-,x^+] \subset \R$ with boundary conditions
$j(x^-,t)=j(x^+,t)=0$ for all $t\geq 0$. The microscopic conservation law 
for $M^\alpha := \bar{\sigma}^\alpha L$ translates
into $\int_{x^-}^{x^+} \rmd x \ \sigma^\alpha(x,t) = \bar{\sigma}^\alpha$.
Defining $\bar{\sigma}^{n+1}:=0$ and the points 
$x^\alpha:= x^+(1-\bar{\sigma}^\alpha) -x^-\bar{\sigma}^\alpha$ one
finds on $[x^-,x^+]$ the static (weak) solution $\sigma_{stat}^\alpha(x) = 
\Theta(x - x^\alpha)$ which one expects as final state of the evolution
of the shock solutions described above.

This implies for the individual static
densities $0\leq \rho_{stat}^\alpha(x) \leq 1$ 
\bel{statdensityprofile}
\rho_{stat}^\alpha(x) = \Theta(x - x^\alpha)
\Theta(x^{\alpha+1} - x) \\
= \left\{\ba{ll} 1 & x^\alpha < x < x^{\alpha+1} \\ 0 & \mbox{ else.} \ea \right. 
\ee 
This is a phase-separated state with 
successive blocks fully occupied by particles of species $\alpha$ in 
increasing order, reminiscent of phase separation in a related
two-species exclusion processes \cite{Evan98b,Arnd98b,Clin03,Ayer09,Bodi11}.
Also this feature will be further elucidated below.

\section{Results}

\subsection{Reversible measures}

All invariant measures of the $n$-priority ASEP on the finite lattice
$\Lambda^L$ are obtained in explicit form and blocking measures for the
infinite system are presented.

\begin{theo}
\label{Theo:revmeas}
(Reversible measures) 
For $\Lambda=\Lambda_{L}$ and for $\bfeta\in\S_{0,n}^L$ let
\bel{Def:Energy}
E(\bfeta) := - \sum_{k=L^-}^{L^+} \sum_{l=L^-}^{k-1} \sgn(\eta_k-\eta_l)
\ee
be the ``energy'' of a configuration $\bfeta$. Then for
\bel{revmeas}
\pi(\bfeta) :=  q^{-E(\bfeta)}.
\ee
we have:\\

\noindent (i) $\pi$ is a reversible measure for the $n$-priority ASEP \eref{generator} 
on $\Lambda_{L}$.\\

\noindent (ii) With particle numbers $\vec{N} =\{N^1,\dots ,N^n\}$ for each species $\alpha\geq 1$ the canonical measure
\bel{canmeas}
\pi^c_{\vec{N}}(\bfeta) = \pi(\bfeta)
\frac{\mathcal{P}_{\vec{N}}(\bfeta) }{C_{L}(\vec{N})} 
\ee
with the canonical partition function
\bel{canpartfun}
C_{L}(\vec{N}) = \frac{[L]_q!}{\prod_{\alpha=0}^n [N^\alpha]_q!}
\ee
is the unique invariant measure for the process \eref{generator} on the subset 
$\S^L_{\vec{N}}$.\\

\noindent (iii) With chemical potentials $\mu_\alpha \in \R$ and 
$\vec{\mu} := (\mu_1,\dots,\mu_n)$ the grand-canonical family
\bel{grandmeas}
\pi^g_{\vec{\mu}}(\bfeta)= 
\frac{\rme^{\sum_{\alpha=1}^n \mu_\alpha N^\alpha(\bfeta)}}{Z_{L}(\vec{\mu})} 
 \pi(\bfeta)
\ee
are invariant measures with the homogeneous multivariate 
Rogers-Szeg\H{o} polynomial \cite{Gasp04}
\bel{grandpartfun}
Z_{L}(\vec{\mu}) = \sum_{\vec{N}} \rme^{\sum_{\alpha=1}^n \mu_\alpha N^\alpha} 
C_{L}(\vec{N})
\ee
as grand-canonical partition function.
\end{theo}

\begin{rem}
For $n=1$ one has $Z_{L}(\mu)= \prod_{k=L^-}^{L^+}(1+\rme^\mu q^{2k-L^+-L^-})$ and
the grandcanonical invariant measure \eref{grandmeas} becomes a product measure
with marginal densities 
$\rho_k = \rme^\mu q^{2k-L^+-L^-}/(1+\rme^\mu q^{2k-L^+-L^-})$
which is the blocking measure of \cite{Ligg85}
for the ASEP restricted to the subset $\Lambda_{L}$ of $\Z$. The parameter
$\mu$ fixes the center of the shock, which is the lattice point
$k^\ast := \min_{k}{(\rho_k>1/2)}$.
\end{rem}

We point out alternative forms of writing the reversible measure.
From \eref{sigma} we define the partial energies
\be 
\label{Def:Eab} 
E^{\alpha\beta}(\bfeta)  :=  - \sum_{k=L^-}^{L^+} \sum_{l=L^-}^{k-1}
(n_k^\alpha(\bfeta) n^\beta_l(\bfeta) - n_l^\alpha(\bfeta) n^\beta_k(\bfeta)) .
\ee
By sum rules \eref{sumrule3} and \eref{resum2} we have
\bel{Eab}
E^{\alpha\beta}(\bfeta) 
= - \sum_{k=L^-}^{L^+} n_k^\alpha(\bfeta) N_k^\beta(\bfeta) 
= \sum_{k=L^-}^{L^+} n_k^\beta(\bfeta) N_k^\alpha(\bfeta).
\ee
We decompose the energy as
\bel{Edec}
E(\bfeta) = \sum_{\alpha=1}^n \sum_{\beta=0}^{\alpha-1} E^{\alpha\beta}(\bfeta) 
= E^0(\bfeta) + \bar{E}(\bfeta)
\ee
where, both in occupation variable and coordinate representation,
\bea
\label{E0}
E^0(\bfeta) & := & - \sum_{k=L^-}^{L^+} N^0_k(\bfeta) 
= - \sum_{i=1}^{N^0(\bfeta)}(L^+ + L^- - 2x^0_i(\bfeta))\\
\label{E0co}
E^0(\bfx) & = & - \sum_{i=1}^{N(\bfx)}(2x_i - L^+ - L^-)
\eea
and
\bea  
\label{Ered}
\bar{E}(\bfeta) & := & \sum_{\alpha=1}^n \sum_{\beta=1}^{\alpha-1} 
E^{\alpha\beta}(\bfeta) \\
\label{Eredco}
\bar{E}(\bfx)  & = & - \sum_{i=1}^{N(\bfx)} \sum_{j=1}^{i-1}
\sgn(\alpha_i-\alpha_j) .
\eea
This decomposition leads to a factorization of the reversible measure
\bel{revmeasfac}
\pi(\bfeta) = \pi^0(\bfeta) \bar{\pi}(\bfeta) = \pi^0(\vec{x})
\bar{\pi}(\vec{\alpha})
\ee
where the reduced measure
\be 
\label{pired}
\bar{\pi}(\bfeta):= q^{-\bar{E}(\bfeta)}
\ee
does not depend on the vacancy projectors.

Notice that $-L^2 n/(2n+2) \leq E \leq L^2 n/(2n+2)$. For fixed particle numbers
$N^\alpha$ the minimal energy $E^{min}_{\vec{N}} =
- \sum_{\alpha=1}^n\sum_{\beta=0}^{\alpha-1} N^\alpha N^\beta$ 
is achieved when all two-particle signs
are positive, which is the case for the block configuration $\bfeta^{min}$ 
with all vacancies to the left, followed by blocks of species $\alpha$ in increasing
order such that
\bel{etamin}
\bfeta^{min}_k = \alpha \mbox{ for } L^+-M^\alpha < k \leq L^+-M^{\alpha+1}
\ee
with the convention $M^{n+1}=0$. This is a microscopic realization of
the phase-separated macroscopic stationary
density profile \eref{statdensityprofile}. The local energy associated with site
$k$ is long-ranged even though the stochastic dynamics are local. A similar phenomenon
was found in the ABC-model \cite{Evan98b,Clin03,Ayer09,Bodi11} which also 
exhibits phase separation.

Indeed, the discussion of the 
hydrodynamic limit suggests that on microscopic scale the 
invariant infinite volume measure inside the blocks, i.e., on 
points $k = [x]$ with $x^\alpha  < x < x^{\alpha+1}$ is 
concentrated on the configurations with all sites
occupied by a particle of species $\alpha$. In the (lattice) vicinity of the phase boundaries $k^\alpha := [x^\alpha]$
one expects the invariant infinite-volume measure to be blocking measures 
with particles of species $\alpha$ and $\alpha+1$.
This intuition is well borne out by
the form of the blocking measure for the standard ASEP ($n=1$) \cite{Ligg85}
and by the fact that the $n$-species priority ASEP evolving on a 
subset of $\S_{0,n}^{L}$ with only two species $\alpha, \beta \in \S_{0,n}$ 
of any kind an be identified
with the standard ASEP, see Remark \eref{redproc}. 
Thus one obtains infinite-volume blocking
measures of the $n$-species priority ASEP:

\begin{theo}
\label{Theo:Blocking} (Blocking measures)
Let $\alpha, \beta \in S_{n+1}$ be two different particles species of any kind with
$\alpha < \beta$ and fix a strictly positive constant $\lambda \in \R_+$. Then
the product measures $\pi_{\lambda}^{\alpha\beta}$ 
with marginals
\bel{blockmeas}
\rho^\gamma_k := \Prob{n^\gamma_k = 1} = \left\{ \ba{ll}
\frac{1}{1 + \lambda q^{2k}} & \mbox{ for } \gamma=\alpha \\[2mm]
\frac{\lambda q^{2k}}{1 + \lambda q^{2k}} & \mbox{ for } \gamma=\beta \\[2mm]
0 & \mbox{ else } \ea \right.
\ee
are invariant measures for the $n$-species priority ASEP on $\Z$.
\end{theo}

\begin{rem}
The particle density $\rho^\beta_k$ of species $\beta$
has the shape of a discretized shifted hyperbolic tangent
with $\lim_{k\to-\infty} \rho^\beta_k = 0$ and $\lim_{k\to\infty} \rho^\beta_k = 1$.
The parameter $\lambda$ determines the lattice point $k$ where the local particle
density $\rho^\beta_k$ comes closest to 1/2. Specifically, for $\lambda=q^{-2k_0}$
one has $\rho^\alpha_{k_0} = \rho^\beta_{k_0} = 1/2$.
\end{rem}

\subsection{Duality}

We establish that the $n$-priority ASEP is self-dual w.r.t. a family of
duality functions which arise from the symmetry of the generator of the
process under the quantum algebra $U_q[\mathfrak{gl}(n+1)]$.

\begin{theo}
\label{Theo:duality} (Self-duality)
Fix arbitrary parameters $c_\alpha \in \C$, $\alpha\in\S_{0,n}$ and
let $\bfeta$ and $\bzeta$ be two configurations of the $n$-component priority 
ASEP defined by \eref{generator} on the finite lattice $\Lambda_{L}$.
The process is self-dual with respect to the duality function
\bel{theo2}
D(\bzeta,\bfeta;\vec{c}) =  
\prod_{k=L^-}^{L^+} \prod_{\alpha=0}^{n}  (Q^\alpha_{k}(\bfeta)q^{
c_\alpha(M^\alpha(\bfeta)-L)})^{n^\alpha_k(\bzeta)}  
\ee
where 
\bel{Qi}
Q^\alpha_k(\bfeta) = m_k^\alpha(\bfeta) q^{M_k^\alpha(\bfeta)} 
\ee
with the indicator functions $n_k^\alpha(\cdot)$, 
$m_k^\alpha(\cdot)$ \eref{Def:occind} and the 
particle balance $M_k^\alpha(\cdot)$ \eref{Def:Nalphafun}.
\end{theo}

\begin{rem}
\label{Rem:duality}
The duality function $D(\bzeta,\bfeta;\vec{c})$ is neither the duality 
function of \cite{Beli15b} for $n=2$ nor the duality function of 
\cite{Kuan16,Kuan16b} but a new
duality function for the ASEP with second-class particles.
A complete classification of duality functions is not the purpose of this work. 
However, the algebraic methods employed in the proof (which is constructive)
exhibit quite explicitly how other duality functions can be computed, 
including the duality functions of \cite{Beli15b} and \cite{Kuan16,Kuan16b}, 
see Remark \eref{Rem:Uniqueness} for further details.
\end{rem}

We point out alternative forms of writing the duality function and the
corresponding duality matrix. They follow
from \eref{theo2} by using
$m^0_k=1$ which gives $Q^0_k(\bfeta) = q^{2k-L^+-L^-}$, the sum rules \eref{sumrule1}
and \eref{resum1}
which yield 
\be 
\prod_{k=L^-}^{L^+}q^{2k-L^+-L^-}=1, \quad
\prod_{k=L^-}^{L^+}q^{(2k-L^+-L^-)n^0_k}=\prod_{k=L^-}^{L^+}q^{-N^0_k},
\ee 
and the projector property of Lemma \eref{projlem} which implies 
$m_k^\alpha q^{m_k^\alpha} = q m_k^\alpha$. Thus one has
\bea
\label{theo2a}
D(\bzeta,\bfeta;\vec{c}) & = & \prod_{k=L^-}^{L^+}  q^{-N^0_k(\bzeta)}
\prod_{k=L^-}^{L^+} \prod_{\alpha=1}^{n}  (Q^\alpha_{k}(\bfeta)q^{
c_\alpha(M^\alpha(\bfeta)-L)})^{n^\alpha_k(\bzeta)}  \\
\label{theo2b}
& = & \prod_{k=L^-}^{L^+} \prod_{\alpha=1}^{n}  
(\tilde{Q}^\alpha_{k}(\bfeta;c_\alpha))^{n^\alpha_k(\bzeta)}
\eea
where 
\bel{tQi}
\tilde{Q}^\alpha_{k}(\bfeta;c_\alpha) = 
q^{-(1+c_\alpha)\sum_{l=L^-}^{k-1}(1-m_l^\alpha(\bfeta)) + 
(1-c_\alpha)\sum_{l=k+1}^{L^+}(1-m_l^\alpha(\bfeta))}
m_k^\alpha(\bfeta).
\ee
Specifically when $c_\alpha=c_n$ for all $\alpha$ and
and $c_n$ takes values $0$ or $\pm 1$ we denote
$\tilde{Q}^\alpha_{k}(\bfeta) :=  \tilde{Q}^\alpha_{k}(\bfeta;0)$,
$D(\bzeta,\bfeta) := D(\bzeta,\bfeta;\vec{0})$,
$D^\pm(\bzeta,\bfeta)  :=  D(\bzeta,\bfeta;\pm \vec{1})$.

In coordinate representation one has
for a configuration $\bfx \in V^L_{\vec{N}}$ with $N(\bfx)$ particles
\be 
\label{theo2bc}
D(\bfx,\bfeta;\vec{c})
= \prod_{i=1}^{N(\bfx)}
\tilde{Q}^{\alpha_i}_{x_i}(\bfeta;c_{\alpha_i}) =: Q^{\vec{c}}_{\bfx}(\bfeta)
\ee
with the particle coordinates $x_i$ and species
$\alpha_i$ and therefore
\bea
\label{D}
D(\bfx,\bfeta)
& = & \prod_{i=1}^{N(\bfx)}
\tilde{Q}^{\alpha_i}_{x_i}(\bfeta)  =: Q_{\bfx}(\bfeta)\\
\label{theo2a1}
D^+(\bfx,\bfeta)
& = & \prod_{i=1}^{N(\bfx)} 
q^{- 2 \sum_{l=L^-}^{x_i-1} (1-m_l^{\alpha_i}(\bfeta))}
m^{\alpha_i}_{x_i}(\bfeta)  =: Q^+_{\bfx}(\bfeta) \\
\label{theo2a2}
D^-(\bfx,\bfeta)
& = & \prod_{i=1}^{N(\bfx)} 
q^{2 \sum_{l= x_i+1}^{L^+} (1-m_l^{\alpha_i}(\bfeta))}
m^{\alpha_i}_{x_i}(\bfeta)  =: Q^-_{\bfx}(\bfeta).
\eea

Notice that in the duality functions $D^\pm$ 
the lattice enters through
the position of only one its edges $L^\pm$ in the exponential factors.
This observation 
together with Lemma \eref{Lemma:ASEPZ} and the monotonicity property
$1-m_l^\alpha(\bfeta) \leq 1-n_l^n(\bfeta)$ for all $\alpha$ and all $l$ 
yields duality
functions for the semi-infinite lattices $\Lambda_{\pm}$ and the infinite
integer lattice $\Z$:

\begin{cor}
Let $\bfx$ be the coordinate representation of a configuration
with a finite number $N$ 
of particles of each species $\alpha \geq 1$ and let 
$\Lambda$ be a lattice defined in \eref{Def:lattice}. Then:\\

\noindent a) The $n$-priority ASEP on $\Lambda_{\pm}$  
is self-dual w.r.t. the duality function $D^\pm(\bfx,\bfeta)$ 
(\eref{theo2a1} or \eref{theo2a2} resp.) for all $\bfeta
\in \S_{0,n}^{\pm}$.\\

\noindent b) The $n$-priority ASEP on $\Z$  
is self-dual w.r.t. the duality function $D^+(\bfx,\bfeta)$ \eref{theo2a1} 
($D^-(\bfx,\bfeta)$ \eref{theo2a2}) for all $\bfeta
\in \S^{\Z}_{n^{>}}$ \eref{rfull} ($\S^{\Z}_{n^{<}}$ \eref{lfull}).
\end{cor}

Interpreting the duality function 
$Q^{\vec{c}}_\bfx(\bfeta)$ as a function 
of the configurations $\bfeta$ indexed by $\bfx$ (without 
thinking of $\bfx$ as representing another particle configuration)
the self-duality relation \eref{theo2} reads 
\bel{theo2c}
\E^{\bfeta} Q^{\vec{c}}_\bfx(\bfeta_t) = \sum_{\bfy \in V^\Lambda_{\vec{N}(\bfx)}}
Q^{\vec{c}}_\bfy(\bfeta) P(\bfy,t|\bfx,0).
\ee
with the
transition probability $P(\bfy,t|\bfx,0)$ from configuration $\bfx$ at time $t=0$
to $\bfy$ at time $t\geq 0$.
This transition probability 
can be computed explicitly from the nested coordinate Bethe ansatz 
\cite{Yang67,Gaud14} for the $n$-component 
priority ASEP, see \cite{Dahm95,Popk02,Cant08,Chat10} for work on $n=2$ 
and \cite{Arit09,Trac13} for general $n$.
Specifically for $N(\bfx)=1$ this observation leads to

\begin{cor}
Let $\bfx$ be the configuration with a single particle of some species 
$\alpha \geq 1$ at 
site $k$, $\bfeta\in\S_{0,n}^L$. Then 
\bel{continuity3}
\mathcal{L} Q^\alpha_k(\bfeta) = 
J^\alpha_{k-1}(\bfeta) - J^\alpha_k(\bfeta) , \quad L^- \leq k \leq L^+ 
\ee
with the locally conserved instantaneous currents
\bel{conscurr3} 
J^\alpha_k(\bfeta) = w \left(q Q^\alpha_{k}(\bfeta) 
- q^{-1} Q^\alpha_{k+1}(\bfeta)\right).
\ee
for $L^- \leq k \leq L^+$ and $J^\alpha_{L^--1}=J^\alpha_{L^++1}=0$ for all
$\alpha \in \S_{0,n}$.
\end{cor}

%\begin{rem} 
The linearity of the locally conserved currents \eref{conscurr3}
corresponds to a discrete biased diffusion equation for the expectation of 
$Q_k^\alpha$ for arbitrary initial configurations $\bfeta$. One can interpret 
the substitution $m_k^\alpha(\bfeta) \to Q_k^\alpha(\bfeta)$ as the lattice 
analogue of the Cole-Hopf transformation for the Kardar-Parisi-Zhang equation 
\cite{Kard86,Sasa10,Amir11}, thus generalizing an analogous earlier observation 
for $n=1$ \cite{Bert97,Schu97}.
%\end{rem}

\subsection{Microscopic structure of shocks}

Here we show, using self-duality, that on microscopic scale shocks
in the $n$-priority ASEP on $\Z$, or more precisely microscopic shock markers 
that indicate the positions of shocks on macroscopic scale, perform a 
random motion that we call the {\it shock exclusion process}.

The shock exclusion process, defined formally below, has jumps like
the $n$-priority ASEP, but each particle
has its individual jump rate. Thus the process has a natural description
in terms of coordinate representation: If site $x_i\pm 1$ next
to particle $i$ is empty, then particle $i$ jumps
after an exponentially distributed random time with parameter
$w v^{\pm 1}_i$ defined below in \eref{bareshockhoppingrates}.
If two particles $i,i+1$ are nearest neighbours they
exchange their colour after an exponential random time 
with parameter (a) $wq^{-1}$ if 
$\alpha_i > \alpha_{i+1}>0$, (b) $q$ if $0< \alpha_i < \alpha_{i+1}$, 
(c) $wq$ if $\alpha_{i+1}=0$, and (d) $q^{-1}$ if 
$\alpha_i=0$. Notice that for
colours $\alpha \in \{1,\dots,n-1\}$ the hopping bias $q$ is inverted
compared to the $n$-priority ASEP.

\begin{df} 
\label{Def:shockprocess} (Shock exclusion process)
Let $\bfx =(\vec{x},\vec{\alpha}) \in V^\Z_N $ be a configuration with $N$ particles and 
let for a parameter $\rho_0 \in (0,1]$
the functions $v_i$, $1 \leq i \leq N$, be defined by
\bel{bareshockhoppingrates}
v_i := \frac{(q-q^{-1})\rho_i(1-\rho_i)}{\rho_{i}-\rho_{i-1}}, \quad 
v^{-1}_i = \frac{(q-q^{-1})\rho_{i-1}(1-\rho_{i-1})}{\rho_{i}-\rho_{i-1}}.
\ee

The shock exclusion process is defined by the generator
\bel{shockgen}
\mathcal{M} f(\bfx) = \sum_{i=1}^N \mathcal{D}_i f(\bfx)
+ \sum_{i=1}^{N-1} \mathcal{C}_i f(\bfx)
\ee
with the single-particle hopping and the two-particle colour exchange generators
\be 
\label{shockgenlocal}
\mathcal{D}_i f(\bfx) = \sum_{\sigma = \pm}
w_i^{\sigma}(\bfx)(f(\bfx_i^{\sigma})-f(\bfx)) , \quad
\mathcal{C}_i f(\bfx) = c^{ii+1}(\bfx) (f(\bfx^{ii+1})-f(\bfx)) 
\ee
with the jump rates 
\bel{shockhoppingrates}
w_i^{\pm}(\bfx) := 
w v^{\pm 1}_i(1-\delta_{x_{i\pm 1},x_i\pm 1}),\quad
1\leq i \leq N,
\ee
and for $1\leq i \leq N-1$ the colour exchange rates
\bea
c^{ii+1}(\bfx) & := & \left[ 
w \left( q^{\sgn(\alpha_{i+1}-\alpha_i)} - \delta_{\alpha_i,\alpha_{i+1}} \right)
(1-\delta_{\alpha_i,0})(1-\delta_{\alpha_{i+1},0}) \right.
\nonumber \\
\label{colourrates}
& & + \left. w q (1-\delta_{\alpha_i,0})\delta_{\alpha_{i+1},0}
+w q^{-1} \delta_{\alpha_i,0}(1-\delta_{\alpha_{i+1},0})\right] 
\delta_{x_{i+1},x_i+1}
\eea
where $\bfx_i^{\sigma}$ and $\bfx^{ii+1}$ are defined by
\be 
\ba{ll}
 (\vec{x}_i^{\,\pm})_j
= x_j  \pm \delta_{j,i} & \vec{\alpha}^\pm = \vec{\alpha} \\[4mm]
(\vec{x}^{\,ii+1})_j = x_j & \vec{\alpha}^{\,ii+1}
= \pi^{ii+1}(\vec{\alpha}) 
\ea
\ee
with the conventions $x_0 := -\infty$, $x_{N+1} := +\infty$, and 
the colour permutation $\pi^{ii+1}(\cdot)$ defined analogously to 
\eref{perm}.
\end{df}

\begin{rem}
The definition of $v_i$ implies that 
\bel{shockfugacityratio} 
\frac{\rho_{i}(1-\rho_{i-1})}{\rho_{i-1}(1-\rho_{i})} = q^2, 
\quad  \quad \forall i \in \{1,\dots,N\}.
\ee
and also $w^+_i w^-_i = w^2$ as one has for the single-particle hopping rates
$w^\pm = wq^{\pm 1}$ of the $n$-priority ASEP. In fact, one can show (see Proposition \eref{Prop:Htrafo} below) that the shock exclusion process and
the $n$-priority ASEP on the finite lattice $\Lambda_L$ are related
by a similarity transformation, up to a boundary term.
\end{rem}

For $N=1$ the shock exclusion process reduces to a simple random walk of a 
particle that moves with average speed 
\bea
\label{shockvel1}
v := \E^0\, \frac{x(t)}{t}
& = & w (q-q^{-1}) (1-\rho_1-\rho_0)
\eea
and diffusion coefficient
\bel{shockdiff1}
D := \E^0\, \frac{(x(t)-vt)^2}{t} = \frac{w}{2}(q-q^{-1})
\frac{\rho_1(1-\rho_1)+\rho_{0}(1-\rho_{0})}{\rho_{1}-\rho_{0}}
\ee
and Gaussian fluctuations on coarse-grained diffusive scale.
For $N$ particles the shock exclusion process can be seen on large scales 
as a gas of particles with different masses (or friction coefficients) 
that drift diffusively 
in a gravitational field with mean velocities $v_i > v_{i+1}$ until 
they collide, with $v_i$ and diffusion coefficients $D_i$ given by
\bea
\label{shockvels1}
v_i & = & \frac{w(q-q^{-1})\rho_i(1-\rho_i)}{\rho_{i}-\rho_{i-1}}-
\frac{w(q-q^{-1})\rho_{i-1}(1-\rho_{i-1})}{\rho_{i}-\rho_{i-1}} \\
\label{shockvels2}
& = &  w (q-q^{-1}) (1-\rho_{i}-\rho_{i-1}) \\
\label{shockdiff}
D_i & = & \frac{w}{2}(q-q^{-1})
\frac{\rho_i(1-\rho_i)+\rho_{i-1}(1-\rho_{i-1})}{\rho_{i}-\rho_{i-1}} .
\eea

\begin{df}
\label{Def:shockmeas} (Shock measure)
Let $\bfx \in V^{\Z}_K$ be a configuration 
with $K$ particles at positions
$x_j\in\Z$, $x_{j+1} > x_j$ for all $j \in \{1,\dots,K\}$, with 
$\alpha(j) = \eta_{x_j}(\bfx) \in \S_n$
specifying the particle species at position $x_j$ and define
\bel{shockmargfunction}
\rho_j(\lambda) := \frac{q^{2j-K +\lambda}}{1+q^{2j-K +\lambda}}, 
\quad 0 \leq j \leq K
\ee
with parameter $\lambda\in\R$.
The product measures $\nu_{\bfx}(\cdot)$ on $\S^\Z_{0,n}$ indexed by $\bfx$
%\bel{shockmeasinf}
%\nu_{\bfx} = \prod_{k} \nu_{\bfx}^k
%\ee
with parameter $\lambda$ and marginals $\nu_{\bfx}^k(\eta_k)$ for $\eta_k\in\S_{0,n}$ at site $k\in\Z$ given by
\bel{shockmarg} 
\nu_{\bfx}^k(\eta_k) = \left\{ \ba{ll}
\delta_{\eta_k,\alpha(j)-1} & k = x_j \\
\rho_0(\lambda) \delta_{\eta_k,n} + (1-\rho_0(\lambda))\delta_{\eta_k,0}  & k < x_{1} \\
\rho_j(\lambda) \delta_{\eta_k,n} + (1-\rho_j(\lambda))\delta_{\eta_k,0}  & x_j < k < x_{j+1}, 1 \leq j \leq N-1 \\
\rho_N(\lambda) \delta_{\eta_k,n} + (1-\rho_K(\lambda))\delta_{\eta_k,0}  & k > x_{K}
\ea\right. 
\ee
are called shock measures
for the $n$-priority ASEP on $\Z$ 
with $N^\alpha$ microscopic shock markers of type 
$\alpha(j)$ at positions $x_j$ and $K=\sum_{\alpha=0}^{n-1} N^\alpha$ shocks.
The restriction to $\Lambda_L$ for $\bfx\in V^L_{\vec{N}}$
\bel{shockmeas}
\mu^L_{\bfx}(\bfeta) := \prod_{k=L^-}^{L^+} \nu_{\bfx}^k(\eta_k)
\ee
is also called shock measure.
\end{df}

We can identify the positions $\vec{x}$ and types $\vec{\alpha}$ 
of the shockmarkers in \eref{shockmarg}  
with a configuration $\bfx = (\vec{x},\vec{\alpha^\Gamma})
\in V^\Z_N$, and similarly with $\mu^L_{\bfx}$
for the finite lattice $\Lambda_L$. We shall call this identification
canonical.

%With the definition of the shock exclusion process and the shock measures 
Now we are in a position to state the main result of this subsection.

\begin{theo}
\label{Theo:shock}
Let $\nu_\bfx(t)$ denote
the distribution at time $t$ of the $n$-priority ASEP, starting from an
$N$-particle shock measure $\nu_\bfx$. Then, for any $x\in\Z$ 
\bel{shockevolution}
\nu_{\bfx}(t) = \sum_{\bfy \in V^\Z_N} 
P(\bfy,t|\bfx,0) \, \nu_{\bfy}
\ee
where $P(\bfy,t|\bfx,0)$ is the transition probability 
of the shock exclusion process \eref{Def:shockprocess}
with the canonical identification of the shockmarkers in 
$\nu_{\bfx}$ with configurations in $\bfx\in V^\Z_N$.
\end{theo}

\begin{cor}
For $N$ shock markers the invariant distribution of the distances 
$r_i := x_{i+1} - x_i$ is a product distribution with geometric marginals
\bel{shockdist}
\Prob{(r_i=k)} = \frac{(\rho_N-\rho_i)(\rho_i-\rho_0)}{\rho_i(1-\rho_i)}
 \left(1 - \frac{(\rho_N-\rho_i)(\rho_i-\rho_0)}{\rho_i(1-\rho_i)}\right)^k
\ee
We call this stationary limit a bound state of $N$ shocks, which on
macroscopic scale is a single shock with a jump discontinuity in the 
density from
$\rho_0$ to $\rho_N$.
\end{cor}

This follows by straightforward computation from the standard mapping of the ASEP 
to the zero-range process where the distance process $r_i(t) := x_{i+1}(t) -
x_i(t)$ between $N-1$ neighbouring 
particles in the shock exclusion process is
a zero-range process on $N-1$ sites with site-dependent hopping rates 
\cite{Benj96}. 
With $z=\rho_0/(1-\rho_0)$ we can write
\bel{shockdist2}
\Prob{(r_i=k)} = 1 - \frac{(q^{2N} -q^{2k})(q^{2k}-
1)}{q^{2k}(1+z)(1+q^{2N}z)}.
\ee
With $z = q^{-N-2p}$ the $[(N/2+p)]^{th}$ microscopic 
shock marks the center of 
the macroscopic shock, i.e. $N/2+p = \min_k{[\rho_k>1/2]}$.

Several remarks are in place.

\begin{rem}

\noindent 1) The proof of the theorem uses a specific
duality function. Similar statements for other shock measures
can be obtained for other choices of duality
functions. It is not the purpose of this work to explore this fact 
in more detail.\\

\noindent 2) One recognizes in the microscopic shock velocities 
\eref{shockvels1} the Rankine-Hugoniot speeds \eref{RH}
since $j_i:=w(q-q^{-1}) \rho_i(1-\rho_i)$ is the expectation of the 
particle current to the right of shock $i$ and $j_{i-1}$ is the current 
to the left of shock $i$ in the ASEP. As seen from \eref{shockvels2} 
they coincide with the shock velocity \eref{shockvel} for macroscopic shocks 
of type $\alpha=n$.
The shock diffusion coefficients \eref{shockdiff}
are consistent with the general result
of \cite{Ferr94a} on shock motion in the ASEP on diffusive scale, see
also \cite{Beli02} for $n=1$ and \cite{Bala10} for $n=2$.
The fact that the shocks can be marked with particles of arbitrary
colour of lower priority and the ensuing colour exchange process 
between shock markers is a new insight as well as the point that
the shock exclusion process has its origin in duality.\\

\noindent 3) The transition probability of a generalized shock exclusion process
for $n=1$ and any $N$ has been calculated from Bethe ansatz
in the limit $w\to 0$ and $q\to\infty$ such that $wq=1$ for arbitrary
limiting rates $v_i$ and can be expressed in determinantal form 
\cite{Rako06}.\\

\noindent 4) On macroscopic scale the final stage of the time evolution
of the shock measure can be interpreted as coalescence of shocks
\cite{Ferr00}. In the associated (non-Hermitian) quantum problem this invariant state 
is a many-body bound state of (coloured) magnons in the XXZ-quantum
chain ($n=1$) or their higher rank analogues in the Perk-Schultz chain resp. 
($n>1$). This observation suggests that many-body bound states of other 
quantum systems, that are known to exist from Bethe-ansatz, might have a 
classical analog as shock coalescence or related phase separation
phenomena.
\end{rem}

\section{Proofs}

\subsection{Proof of Theorem \eref{Theo:revmeas}}

(i) Reversibility is guaranteed by the detailed balance condition \eref{detbal} 
which here reads $w^{rr+1}_q(\bfeta^{rr+1}) 
= w^{rr+1}_q(\bfeta) q^{-E(\bfeta)+E(\bfeta^{rr+1})}$
for all $r$ with the total energy \eref{Def:Energy}. 
We observe that
since $\bfeta$ and $\bfeta^{rr+1}$ differ only at sites $r$ and $r+1$ 
by a permutation one has
$E(\bfeta^{rr+1}) - E(\bfeta) =  2\, \sgn(\eta_{r+1}-\eta_r)$.
The detailed balance relation then follows from 
$w^{rr+1}_q(\bfeta^{rr+1}) = w (q^{\sgn(\eta_{r+1}-\eta_r}
- \delta_{\eta_{r+1},\eta_r})$ and $w^{rr+1}_q(\bfeta) = w 
(q^{-\sgn(\eta_{r+1}-\eta_r} - \delta_{\eta_{r+1},\eta_r})$.\\

\noindent (ii) From (i) and particle number conservation it follows that the 
canonical measure \eref{canmeas} is an invariant measure for the process 
$\bfeta_t$. Uniqueness follows from ergodicity for fixed particle numbers which 
itself a consequence of the fact that the process is a sequence of permutations. 
It remains to prove normalization $\sum_{\bfeta} \pi^c_{{N}}(\bfeta) =1$. 

We work with the coordinate representation and 
recall the decomposition \eref{Edec} of the total energy and the corresponding
decomposition \eref{revmeasfac} of the invariant measure. For fixed particle
numbers $N^\alpha$ for all species $\alpha \geq 1$ (and hence also for 
the vacancies $\alpha=0$ with $N^0 = L - M^1$) we first compute, using
\eref{E0} and $M^1=N$,
\be 
\sum_{\vec{x}\in W^L_{N}} \pi(\bfx) = \sum_{\vec{x}^0\in W^L_{N^0}}
q^{-\sum_{i=1}^{N^0} (2x_i^0 - L^+-L^-)} \bar{\pi}(\vec{\alpha}) =
{L \choose N^0}_q.
\ee
The second equality, which is trivially identical to the sum over the positions 
of all particles in the Weyl alcove $W^L_{N}$, was obtained in \cite{Sand94}, 
but actually goes back to a classical result from the theory of partitions \cite{Andr76}. 

Now we observe that the colour array $\vec{\alpha}$ is a configuration of 
a $(n-1)$-priority ASEP on a finite lattice with $L^-=1$ and $L^+=N=L-N^0$ and by 
\eref{revmeas} $\bar{\pi}(\vec{\alpha})$ is a reversible measure for this process.
Thus we can iterate the decomposition \eref{Edec} with particles of species
1 playing the role of the vacancies. This yields
\be 
\sum_{\vec{\alpha}^1\in W^{L-N^0}_{N^1}} \pi(\vec{\alpha}) 
= \sum_{\vec{\alpha}^1\in W^{L-N^0}_{N^1}}
q^{-\sum_{i=1}^{N^1} (2x_i^1 - (L-N^0)-1)} \bar{\pi}(\vec{\alpha}) =
{L-N^0 \choose N^1}_q.
\ee
By further iteration one gets
\be 
\sum_{\bfeta \in \S^L_{\vec{N}}} \pi(\bfeta) =
\prod_{\alpha=0}^{n-1}  {M^\alpha \choose N^\alpha}_q 
\ee
which is equal to the normalization factor \eref{canpartfun}.\\

\noindent (iii) With (ii) we have
$\exp{\left(\sum_{\alpha=1}^n \mu_\alpha N^\alpha\right)} 
C_{L}(\vec{N}) \pi^c_{{N}}(\bfeta)
= \exp{\left[\sum_{\alpha=1}^n \mu_\alpha N^\alpha(\bfeta)\right]} \pi(\bfeta)
\mathcal{P}_{\vec{N}}(\bfeta)$
and therefore
\be
\sum_{\vec{N}} \rme^{\sum_{\alpha=1}^n \mu_\alpha N^\alpha} C_{L}(\vec{N}) \pi^c_{\vec{N}}(\bfeta)
= \rme^{\sum_{\alpha=1}^n \mu_\alpha N^\alpha(\bfeta)} \pi(\bfeta).
\ee
Since $\pi^c_{\vec{N}}$ is normalized taking the sum over all $\bfeta$ yields the
sum \eref{grandpartfun}. \qed

\subsection{Proof of Theorem \eref{Theo:duality}}

The central ingredients in the proof of self-duality 
are the quantum algebra symmetry of the generator
and its reversibility. The concrete form of the duality function
follows from certain matrix elements of the symmetry operators of
$U_q[\mathfrak{gl}(n+1)]$. 
We begin with establishing the quantum algebra symmetry of $H$.

\subsubsection{Quantum algebra symmetry}

\begin{lmm}
\label{Lemma:xtrafo}
Let $\hat{\pi}$ be the matrix form \eref{revmeasmatrix} of the reversible 
measure \eref{revmeas} and fix a parameter $c\in\C$. The representation 
matrices $\sigma^{\alpha,+}_k$ of
$\mathfrak{gl}(n+1)$ acting non-trivially on site $k$ transform under the
diagonal similarity transformation  
\bel{Def:xtrafo} 
\tilde{\sigma}^{\alpha,\pm}_k(c) 
:= \hat{\pi}^{-c} \sigma^{\alpha,\pm}_k \hat{\pi}^c, \quad 1\leq \alpha\leq n,
\ee
as follows
\bel{xtrafo} 
\tilde{\sigma}^{\alpha,\pm}_k(c)
= q^{\pm c(\hat{N}^\alpha(k) + \hat{N}^{\alpha-1}(k))} \sigma^{\alpha,\pm}_k.
\ee
\end{lmm}

\Proof 
For the matrix elements $(\tilde{\sigma}^{\alpha,+}_k)_{\bzeta\bfeta}(c) = 
\bra{\bzeta} \tilde{\sigma}^{\alpha,+}_k(c) \ket{\bfeta}$
we obtain 
\be 
(\tilde{\sigma}^{\alpha,+}_k)_{\bzeta\bfeta}(c) 
= \frac{\pi^c(\bfeta)}{\pi^c(\bzeta)}
\bra{\bzeta} \sigma^{\alpha,+}_k \ket{\bfeta} 
= \frac{\pi^c(\bfeta)}{\pi^c(\bzeta)} n_k^\alpha(\bfeta)
\delta_{\bzeta,\bfeta^{k,-}}.
\ee 
Since $\sigma(n\pm k) - \sigma(n) 
= \pm (\delta_{n,0} + \delta_{n,\mp k}) \pm 2 \sum_{l=1}^{k-1} \delta_{n,\mp l}$
we have that
\be 
\sigma(\eta^{k,-}_l-\eta^{k,-}_m)
= \sigma(\eta_l-\eta_m) - \delta_{l,k} \left( \delta_{\eta_{l},\eta_{m}} + 
\delta_{\eta_{l},\eta_{m}+1} \right) 
+ \delta_{m,k} \left( \delta_{\eta_{l},\eta_{m}} +
\delta_{\eta_{l},\eta_{m}-1} \right).
\ee 
Specifically, for $\eta_k = \alpha$ this yields
\be 
\sigma(\eta^{k,-}_l-\eta^{k,-}_m) 
= \sigma(\eta_l-\eta_m) - 
\delta_{l,k} \left( n^\alpha_m(\bfeta) + n^{\alpha-1}_m(\bfeta) \right) 
+ \delta_{m,k} \left( n^\alpha_l(\bfeta) + n^{\alpha-1}_l(\bfeta) \right) .
\ee 
Therefore
\bea 
-E(\bfeta) +E(\bfeta^{k,-}) & = & \sum_{l=L^-+1}^{L^+} \sum_{m=L^-}^{l-1} 
\delta_{l,k} \left( n^\alpha_m(\bfeta) + n^{\alpha-1}_m(\bfeta) \right) 
- \delta_{m,k} \left( n^\alpha_l(\bfeta) + n^{\alpha-1}_l(\bfeta) \right) \nonumber \\
& = & \sum_{m=L^-}^{k-1} \left( n^\alpha_m(\bfeta) + n^{\alpha-1}_m(\bfeta) \right)
- \sum_{m=k+1}^{L^+} \left( n^\alpha_m(\bfeta) + n^{\alpha-1}_m(\bfeta) \right)
\eea
which yields \eref{xtrafo} for positive sign.
Going through similar steps one verifies \eref{xtrafo} for $\tilde{\sigma}^{\alpha,-}_k(c)$.
\qed

\begin{prop}
\label{Lem:qsym}
Let $\hat{\pi}$ be the diagonal matrix form of the reversible measure
\eref{revmeas}.
The generator $H$ \eref{H} of the $n$-species priority ASEP
is related to the Perk-Schultz-Hamiltonian $H^{PS}$ 
\eref{HPS} by the ground state transformation 
\bel{similarity}
H^{PS} = (\hat{\pi})^{-1/2} H (\hat{\pi})^{1/2}
\ee
and satisfies the commutation relations
\bel{qsymmetry}
\comm{H}{\hat{N}^\alpha} = \comm{H}{Y^\pm_\alpha} = 0
\ee
with the representation matrices of $U_q[\mathfrak{gl}(n+1)]$ 
given by $\hat{N}^\alpha$ \eref{coprodrepN}
and $Y^\pm_\alpha$ defined by
\bel{Def:Y}
Y^\pm_\alpha := \sum_{k=L^-}^{L^+} Y^{\alpha,\pm}(k)
\ee
with
\bea
\label{coprodyp}
Y^{\alpha,+}(k) & = & 
q^{-\sum_{l=L^-}^{k-1} \hat{n}^\alpha_l} \sigma^{\alpha,+}_k 
q^{\sum_{l=k+1}^{L^+} \hat{n}^\alpha_l} 
 \\
\label{coprodym}
Y^{\alpha,-}(k) 
& = & q^{\sum_{l=L^-}^{k-1} \hat{n}^{\alpha-1}_l}  \sigma^{\alpha,-}_k 
q^{-\sum_{l=k+1}^{L^+} \hat{n}^{\alpha-1}_l} .
\eea
and the fundamental representation $\hat{n}^\alpha$, $\sigma^{\alpha,\pm}$  \eref{fundgln} 
of $\mathfrak{gl}(n+1)$.
\end{prop}

\begin{rem}
The result implies that the coproduct defined via \eref{coprodrepN} and \eref{Def:Y} 
is an algebra homomorphism for $U_q[\mathfrak{gl}(n+1)]$. Lemma
\eref{Lemma:xtrafo} implies 
$(\hat{\pi}^{-1} Y^\pm_\alpha \hat{\pi})^T = Y^\mp_\alpha$.
\end{rem}

\noindent {\it Proof of Proposition \eref{Lem:qsym}:} (i) Since $H^{PS}$ is symmetric the claim \eref{similarity}
is equivalent to $H^T = \hat{\pi}^{-1} H \hat{\pi}$, which is the
matrix form \eref{detbalH} of reversibility established in Theorem \eref{Theo:revmeas}.\\

\noindent (ii) In order to prove \eref{qsymmetry} we assert that
\bel{symtrafo} 
\hat{N}^\alpha = \hat{\pi}^{1/2} \hat{N}^\alpha \hat{\pi}^{-1/2}, \quad
Y^\pm_\alpha = \hat{\pi}^{1/2} X^\pm_\alpha \hat{\pi}^{-1/2}.
\ee
The first transformation property is trivial since both $\hat{N}^\alpha$
and $\hat{\pi}$ are diagonal matrices.
The second transformation property follows from Lemma \eref{Lemma:xtrafo} with $c=-1/2$
and the definition \eref{coprodrepX} of $X^\pm_\alpha$.
The statement \eref{qsymmetry} then follows from the symmetry \eref{Hsym} of 
$H^{PS}$. \qed

Next we state a result for $n=1$ which is a much simplified proof
of the duality relation of \cite{Schu97}. The main interest in the present context
is the subsequent corollary for general $n$. In order to avoid excessive indexing 
we shall omit the superscript 1 for species 1 
throughout the proposition and its proof on all quantities $Q,n,N$ etc.
related to the configurations $\bfeta, \bzeta \in \S_{0,1}^L$.

\begin{prop} 
\label{n1}
Fix $n=1$ and let $Y_1^+$ be the generator \eref{Def:Y} of $U_q[\mathfrak{gl}(2)]$.
Then 
\bel{S1p}
\Upsilon_1^+ = \sum_{r=0}^{L} \frac{(Y_1^+)^r}{[r]_q!}
\ee
has
matrix elements
\bel{auxn1}
\bra{\bzeta} \Upsilon_1^+ \ket{\bfeta} = \prod_{k=L^-}^{L^+} (Q_{k}(\bfeta))^{n_k(\bzeta)}
\ee 
with $Q_{k}(\bfeta)$ defined in \eref{Qi}.
\end{prop}

\Proof Using the explicit form
\be 
\prod_{k=L^-}^{L^+} (Q_{k}(\bfeta))^{n_k(\bzeta)}
= \prod_{k=L^-}^{L^+} \left[ \left( n_k(\bfeta) \right)^{n_k(\bzeta)}
q^{n_k(\bzeta) N_k(\bfeta)} \right]
\ee
we can write \eref{auxn1} in coordinate representation 
$\bfx=\bzeta$ as
\bel{auxn1a} 
\bra{\bfx} \Upsilon_1^+ \ket{\bfeta} = 
\prod_{i=1}^{N(\bfx)} Q_{x_i}(\bfeta) 
= \prod_{i=1}^{N(\bfx)} q^{N_{x_i}(\bfeta)} n_{x_i}(\bfeta) .
\ee

In order prove this we introduce $\tilde{Y}^{1,+}(k) = q^{\hat{N}_k} \sigma^{1,+}_k$.
From the coproduct \eref{coprodyp}
and the fundamental representation \eref{fundgln} one finds
$\comm{Y^{1,+}(k)}{ \tilde{Y}^{1,+}(l)} = 0$ for all $k,l \in \Lambda_{L}$
which implies $\comm{\Upsilon_1^+}{\tilde{Y}^{1,+}(l)} = 0$ for all 
$l\in \Lambda_{L}$.

We also note that one has $\bra{s} = \bra{0}\Upsilon_1^+$ and
$\bra{s} \sigma^{\alpha,+}_k = \bra{s} \hat{n}^{\alpha}_k$ \cite{Schu97}.
Thus for $\vec{x}\in W^L_{N(\bfx)}$ we can write in terms of ordered products
\bea
& & \bra{s} \prod_{i=1}^{N(\bfx)} \hat{Q}_{x_i} = 
q^{ N(\bfx)-1} \bra{s} \overset{\longrightarrow}{\prod_{i=1}^{N(\bfx)}} 
\tilde{Y}^{1,+}(x_i)
= q^{ N(\bfx)-1} \bra{0} \Upsilon_1^+ \overset{\longrightarrow}{\prod_{i=1}^{N^1(\bfx)}} 
\tilde{Y}^{1,+}(x_i) \nonumber \\
& &
= q^{ N(\bfx)-1} \bra{0} 
\overset{\longrightarrow}{\prod_{i=1}^{N(\bfx)}} \tilde{Y}^{1,+}(x_i) \Upsilon_1^+ 
= \bra{\bfx} \Upsilon_1^+ .
\eea
Taking the inner product with $\ket{\bfeta}$ proves \eref{auxn1a}
since $\hat{Q}_{x_i}\ket{\bfeta}=Q_{x_i}(\bfeta)\ket{\bfeta}$. \qed

Lifting Proposition \eref{n1} to $\C^{d_{n,L}}$ we note the following corollary.
\begin{cor}
\label{auxn} 
Define the matrix sums
\bel{Si}
\Upsilon_\alpha^\pm := \sum_{l=0}^{L} \frac{(Y^\pm_\alpha)^l}{[l]_q!}.
\ee
For $1 \leq \alpha\leq n$ and all $\bzeta,\bfeta \in \S_{0,n}^L$
we have
\be  
\bra{\bxi^r} \Upsilon^+_\alpha \ket{\bxi^s} =
\prod_{k=L^-}^{L^+} 
\bar{A}^\alpha_k(q;\bxi^r,\bxi^s) \bar{B}^\alpha_k(\bxi^r,\bxi^s)
\Delta^\alpha_k(\bxi^r,\bxi^s) 
\ee 
with the functions
\bea
\label{Def:Abaralphak}
& & \bar{A}^\alpha_k(q;\bxi^r,\bxi^s) := q^{n^{\alpha}_k(\bxi^r) N_k^\alpha(\bxi^s)},
\quad
\bar{B}^\alpha_k(\bxi^r,\bxi^s) :=
\left( n^{\alpha}_k(\bxi^s) \right)^{n_k^\alpha(\bxi^r)} \\
\label{Def:Deltaalphak}
& & \Delta^\alpha_k(\bxi^r,\bxi^s) := 
\prod_{\overset{\beta=0}{\beta\neq\alpha-1,\alpha}}^n
\delta_{n^\beta_k(\bxi^r),n^\beta_k(\bxi^s)} .
\eea
\end{cor}

\subsubsection{Symmetry and reversibility}

Self-duality stated in matrix form reads \cite{Sudb95}
\bel{Def:selfduality}
D H = H^T D
\ee
where the duality function $D(\bzeta,\bfeta)$ 
is the matrix element $\bra{\bzeta} D \ket{\bfeta}$. A symmetry of $H$
is defined by the commutation relation
\bel{Def:symmetry}
S H = H S
\ee
with symmetry operator $S$. This, along with reversibility \eref{detbalH}
in matrix form, shows that $\hat{\pi}^{-1} S$ is a 
self-duality matrix, see Theorem 2.6 of \cite{Giar09}.
According to Lemma \eref{qsymmetry} $H$ commutes with the matrix sums \eref{Si}
and with the number operators $N^\alpha$ and products of these matrices 
in arbitrary order. Thus we obtain

\begin{prop}
\label{dualitymatrix}
Define
\be
\Upsilon^+ := \Upsilon^+_1 \Upsilon^+_2 \dots \Upsilon^+_n .
\ee
The matrix 
\bel{D0}
D = \hat{\pi}^{-1} \Upsilon^+
\ee 
is a duality matrix for the $n$-component ASEP defined
by the generator \eref{generator}.
\end{prop}

\begin{rem}
\label{Rem:Uniqueness}
The duality matrix is not unique as any other linear combination of products 
of the symmetry operators gives rise to some duality matrix.
For $n=2$ one can show that the duality matrix in Theorem 3.3 in
\cite{Beli15b} is $\hat{\pi}^{-1} \Upsilon^+_2 
\Upsilon^-_1$ while the duality matrix of Theorem 2.2 in \cite{Kuan16} is 
obtained from $\hat{\pi}^{-1} \Upsilon^-_2 \Upsilon^-_1$.
\end{rem}

Observe that $D^{\vec{c}}(\bzeta,\bfeta)$ \eref{theo2} 
for $c_\alpha \neq 0$ can be decomposed
\be 
D^{\vec{c}}(\bzeta,\bfeta) = \prod_{\alpha=0}^n 
q^{c^\alpha N^\alpha(\bzeta) (M^\alpha(\bfeta)-L)} D(\bzeta,\bfeta)
\ee
with $D(\bzeta,\bfeta)$ given in \eref{D}.
Assume now that 
\bel{aux4}
D(\bzeta,\bfeta) 
= \bra{\bzeta} \hat{\pi}^{-1} \Upsilon^+ \ket{\bfeta}.
\ee
Then for a general choice of parameters $c_\alpha$ the Theorem 
follows directly from Proposition \eref{dualitymatrix} and the preceding discussion
since particle 
number conservation is also a symmetry and one can therefore construct
from \eref{aux4} further duality functions which are arbitrary functions
of the conserved particle numbers $N^\alpha(\bzeta)$ and $N^\alpha(\bfeta)$.
Therefore it remains to prove \eref{aux4}.\\

\subsubsection{Proof of \eref{aux4} for $q=1$}

First we consider the symmetric case $q=1$ where $H=H^{PS}$ commutes with 
tensor representation
$S^{\alpha,\pm} = \sum_{k=L^-}^{L^+} \sigma_k^{\alpha,\pm}$ and $N^\alpha$ 
obtained from the coproduct \eref{coprod} for $q=1$, corresponding to symmetry 
under the Lie algebra $\mathfrak{gl}(n+1)$. We generalize the
approach of \cite{Schu94} for $n=1$ to general $n$.

\begin{lmm}
\label{indlem}
Let $\sigma^{\alpha,+}$ be the fundamental representation matrix
of $\mathfrak{gl}(n+1)$ defined in \eref{fundgln}
and let $(\zeta| \in \C^{n+1}$ ($|\eta)\in \C^{n+1}$) be the canonical row (column)
basis vector with component 1 at position $\zeta$ ($\eta$) and
let $S^{\alpha,+} = \sum_{k=L^-}^{L^+} \sigma_k^{\alpha,+}$ be the tensor representation
obtained from the coproduct \eref{coprod} for $q=1$. Then the following
factorization properties hold:
\bea
\label{Sp}
\rme^{S^{\alpha,+}} & = & \prod_{k=L^-}^{L^+} (1+\sigma_k^{\alpha,+}) \\
\label{mton2}
\bra{\bzeta}%\prod_{k=L^-}^{L^+}
\prod_{\alpha=1}^{\overset{\longrightarrow}{n}} \rme^{S^{\alpha,+}}
%(1+\sigma_k^{\alpha,+})
\ket{\bfeta}
& = & \prod_{k=L^-}^{L^+} \prod_{\alpha=1}^n (m_k^\alpha(\bfeta))^{n_k^\alpha(\bzeta)}.
= \prod_{k=L^-}^{L^+} m_k^{\zeta_k}(\bfeta).
\eea
for the ordered product with the lowest index on the left and
with the indicator functions $n_k^\alpha(\cdot)$ and 
$m_k^\alpha(\cdot)$ \eref{Def:occind}.
\end{lmm}

\Proof The commutation relations
\eref{localcommutator} imply
$\exp{(S^{\alpha,+})} = \prod_{k=L^-}^{L^+} \exp{(\sigma^{\alpha,+}_k)}$.
Observing that $\sigma^{\alpha,+}$ (and hence $\sigma^{\alpha,+}_k$) is nilpotent of
degree 2 leads to the first property \eref{Sp}. From \eref{inprodfac} 
in Lemma \eref{Lem:inprodprops} we have the factorization property 
\be 
\bra{\bzeta}\prod_{k=L^-}^{L^+}
\prod_{\alpha=1}^{\overset{\longrightarrow}{n}}(1+\sigma_k^{\alpha,+})\ket{\bfeta}
= \prod_{k=L^-}^{L^+}
(\zeta_k|\prod_{\alpha=1}^{\overset{\longrightarrow}{n}}(1+\sigma_k^{\alpha,+})|\eta_k)
\ee
Using $n_k^\alpha(\bzeta) = \delta_{\zeta_k,\alpha}$ \eref{occupos}
it is therefore sufficient
to prove the lemma for a single fixed $k$. Dropping the subscript $k$ one finds
from the explicit form $\sigma^{\alpha,+}=|\alpha-1)(\alpha|$ 
\be 
(\zeta|\prod_{\alpha=1}^{\overset{\longrightarrow}{n}}(1+\sigma^{\alpha,+}) 
= \sum_{\alpha=\zeta}^n (\alpha|
\ee
and with the definitions \eref{Def:occind} the inner product with $|\eta)$ gives
\be 
(\zeta|\prod_{\alpha=1}^{\overset{\longrightarrow}{n}}(1+\sigma^{\alpha,+}) |\eta)
= m^{\zeta}(\eta).
\ee
On the other hand, again by definition \eref{Def:occind}
\be 
\prod_{\alpha=1}^n (m^\alpha(\eta))^{n^\alpha(\zeta)}
= \prod_{\alpha=1}^n (m^\alpha(\eta))^{\delta_{\alpha,\zeta}}
\ee
which is trivially equal to $m^{\zeta}(\eta)$. \qed

With this result established we return to Proposition \eref{dualitymatrix}
and note that for $q=1$ one has $\pi(\bfeta)=1$. Thus Theorem
\eref{Theo:duality} reduces to 
\bel{theo2q1}
D^0(\bzeta,\bfeta) =  
\prod_{k=L^-}^{L^+} \prod_{\alpha=0}^{n} 
(m^\alpha_k(\bfeta))^{n^\alpha_k(\bzeta)} =: D_{sym}(\bzeta,\bfeta)
\ee
and we need to prove only 
$D_{sym}(\bzeta,\bfeta) = 
\bra{\bzeta} \Upsilon_1^+ \Upsilon_2^+ \dots \Upsilon_n^+\ket{\bfeta}$
for $q=1$. Since by definition $\Upsilon_\alpha^+ = \rme^{S^+_\alpha}$ for $q=1$
this follows from Lemma \eref{indlem}. \qed

Before proceeding to the general case $q\neq 1$
we note the following alternative representation of the duality function 
based on the
expansion
\bel{dualityexpansion} 
\bra{\bzeta} \Upsilon^+ \ket{\bfeta} 
= \sum_{\bxi^1} \dots \sum_{\bxi^{n-1}}
\bra{\bzeta} \Upsilon_1^+ \ket{\bxi^1} \bra{\bxi^1} \Upsilon_2^+ \ket{\bxi^2} \dots
\bra{\bxi^{n-1}} \Upsilon_n^+ \ket{\bfeta}
\ee
for general $q\neq 0$.
Defining for $1\leq \alpha \leq n$
\bea
\label{Def:Deltak}
\Delta_k(\bxi^0,\dots,\bxi^n) & := & \prod_{\alpha=1}^n
\Delta^\alpha_k(\bxi^{\alpha-1},\bxi^{\alpha}) = \prod_{\alpha=1}^n
\prod_{\overset{\beta=0}{\beta\neq\alpha-1,\alpha}}^n
\delta_{n^\beta_k(\bxi^{\alpha-1}),n^\beta_k(\bxi^{\alpha})}  \\
\label{Def:Bbark}  
\bar{B}_k(\bxi^0,\dots,\bxi^n) & := & 
\prod_{\alpha=1}^n \bar{B}^\alpha_k(\bxi^0,\dots,\bxi^n)
= \prod_{\alpha=1}^n \left( n^{\alpha}_k(\bxi^\alpha) \right)^{n_k^\alpha(\bxi^{\alpha-1})}
%\bar{A}_k(q;\bxi^0,\dots,\bxi^n)\Delta_k(\bxi^0,\dots,\bxi^n)
\eea
we have from Corollary \eref{auxn} for $q=1$
\bel{Spq1} 
D_{sym}(\bzeta,\bfeta)
= \sum_{\bxi^1} \dots \sum_{\bxi^{n-1}}
\prod_{k=L^-}^{L^+} \bar{B}_k(\bxi^0,\dots,\bxi^n)
\Delta_k(\bxi^0,\dots,\bxi^n) .
\ee
with $\bzeta=\bxi^0$ and $\bfeta=\bxi^n$.

\subsubsection{Proof of \eref{aux4} for $q\neq 1$}

\noindent $\bullet$ {\it Step 1:} (Vacancy contribution) 
We treat the $Q^0$-contribution to the duality
function separately, using the factorization \eref{revmeasfac}
of the reversible measure. From \eref{theo2a} and \eref{E0} we have
$D(\bzeta,\bfeta) =  q^{E^0(\bzeta)} \bar{D}(\bzeta,\bfeta)$
with
\be 
\bar{D}(\bzeta,\bfeta) = \prod_{k=L^-}^{L^+} 
\prod_{\alpha=1}^{n}  (Q^\alpha_{k}(\bfeta))^{n^\alpha_k(\bzeta)}  .
\ee
We conclude that \eref{aux4} is equivalent to
\bel{aux5}
\bra{\bzeta} \Upsilon^+ \ket{\bfeta}
= \bar{\pi}(\bzeta) \bar{D}(\bzeta,\bfeta)
\ee
which we set out to prove in the next two steps.\\

\noindent $\bullet$ {\it Step 2:} (Matrix elements of symmetry generators) 
Before dealing with \eref{aux5} we prove a technical lemma involving matrix 
elements of the generators of $U_q[\mathfrak{gl}(n+1)]$.

\begin{lmm}
\label{techlem}
Define for $1\leq \alpha \leq n$
\bea
\label{Def:Abark}
\bar{A}_k(q;\bxi^0,\dots,\bxi^n) & := & \prod_{\alpha=1}^n
\bar{A}^\alpha_k(q;\bxi^{\alpha-1},\bxi^{\alpha}) 
\\
\label{Def:Fbar}
\bar{F}(q;\bxi^0,\dots,\bxi^n) & := & \prod_{k=L^-}^{L^+} \bar{A}_k(q;\bxi^0,\dots,\bxi^n)
\Delta_k(\bxi^0,\dots,\bxi^n).
\eea
Then we have that
\bel{aux6} 
\bar{F}(q;\bxi^0,\dots,\bxi^n) = 
\bar{C}(q;\bxi^0,\bxi^n) \prod_{k=L^-}^{L^+} \Delta_k(\bxi^0,\dots,\bxi^n).
\ee
with
\bel{aux7} 
\bar{C}(q;\bzeta,\bfeta) = 
\bar{\pi}(\bzeta)
\prod_{k=L^-}^{L^+} \prod_{\alpha=1}^n q^{n^{\alpha}_k(\bzeta)M_k^\alpha(\bfeta)} .
\ee
\end{lmm}

\Proof We have to take care of the Kronecker-symbols in $\Delta_k(\cdot)$
\eref{Def:Deltak} to substitute the arguments of the particle number functions 
$\bar{A}_k(\cdot)$ in \eref{Def:Fbar}. We decompose
\be 
\prod_{\overset{\beta=0}{\beta\neq\alpha-1,\alpha}}^n
\delta_{n^\beta_k(\bxi^{\alpha-1}),n^\beta_k(\bxi^{\alpha})} =
\prod_{\beta=0}^{\alpha-2}
\delta_{n^\beta_k(\bxi^{\alpha-1}),n^\beta_k(\bxi^{\alpha})}
\prod_{\beta=\alpha+1}^n
\delta_{n^\beta_k(\bxi^{\alpha-1}),n^\beta_k(\bxi^{\alpha})}.
\ee
Thus
\be 
\prod_{\alpha=1}^{n}\prod_{\beta=0}^{\alpha-2}
\delta_{n^\beta_k(\bxi^{\alpha-1}),n^\beta_k(\bxi^{\alpha})}
= \prod_{\alpha=2}^{n}\prod_{\beta=0}^{\alpha-2}
\delta_{n^\beta_k(\bxi^{\alpha-1}),n^\beta_k(\bxi^{\alpha})}
= \prod_{\beta=0}^{n-2} \prod_{\alpha=\beta}^{n} 
\delta_{n^\beta_k(\bxi^{\alpha-1}),n^\beta_k(\bxi^{\alpha})}
\ee
which implies
\bel{Deltares1} 
n^\beta_k(\bxi^{\alpha}) = n^\beta_k(\bxi^{n}) \mbox{ for } 
0 \leq \beta \leq n-2 \mbox{ and } \beta < \alpha < n
\ee
in the $q$-dependent prefactor $\bar{A}_k(q;\bxi^0,\dots,\bxi^n)$ of \eref{Def:Fbar}.
Similarly we note
\be 
\prod_{\alpha=1}^{n}\prod_{\beta=\alpha+1}^n
\delta_{n^\beta_k(\bxi^{\alpha-1}),n^\beta_k(\bxi^{\alpha})}
= \prod_{\alpha=1}^{n-1}\prod_{\beta=\alpha+1}^{n}
\delta_{n^\beta_k(\bxi^{\alpha-1}),n^\beta_k(\bxi^{\alpha})}
= \prod_{\beta=2}^{n}\prod_{\alpha=0}^{\beta-2}
\delta_{n^\beta_k(\bxi^{\alpha}),n^\beta_k(\bxi^{\alpha+1})}
\ee
which implies 
\bel{Deltares2} 
n^\beta_k(\bxi^{\alpha}) = n^\beta_k(\bxi^{0}) \mbox{ for } 
2 \leq \beta \leq n \mbox{ and } 0 < \alpha < \beta.
\ee

Therefore, defining $m^{n+1}(\cdot)=0$, and using the shorthand
$\bar{F}\equiv\bar{F}(q;\bxi^0,\dots,\bxi^n)$
in the following chain of equations, we have the chains of equations
\bea
\bar{F} & \overset{\eref{Def:Abaralphak}}{=} & 
\prod_{k=L^-}^{L^+}\prod_{\alpha=1}^n 
q^{n^{\alpha}_k(\bxi^{\alpha-1}) N_k^\alpha(\bxi^{\alpha})} 
\Delta_k(\bxi^0,\dots,\bxi^n) \\
& \overset{\eref{Deltares2}}{=} & \prod_{k=L^-}^{L^+}\prod_{\alpha=1}^n 
q^{n^{\alpha}_k(\bxi^{0}) N_k^\alpha(\bxi^{\alpha})} 
\Delta_k(\bxi^0,\dots,\bxi^n) \\
& \overset{\eref{sumrule2}}{=} & \prod_{k=L^-}^{L^+}\prod_{\alpha=1}^n 
q^{-n_k^\alpha(\bxi^{\alpha}) N^{\alpha}_k(\bxi^{0})} 
\Delta_k(\bxi^0,\dots,\bxi^n) \\
& \overset{\eref{Def:occind}}{=} & \prod_{k=L^-}^{L^+}\prod_{\alpha=1}^n 
q^{[\sum_{\beta=0}^{\alpha-1} n_k^\beta(\bxi^{\alpha})-1]
N^{\alpha}_k(\bxi^{0}) + m_k^{\alpha+1}(\bxi^{\alpha}) N^{\alpha}_k(\bxi^{0})} 
\Delta_k(\bxi^0,\dots,\bxi^n) \\
& \overset{\eref{Deltares1}}{=} & \prod_{k=L^-}^{L^+}\prod_{\alpha=1}^n 
q^{[\sum_{\beta=0}^{\alpha-1} n_k^\beta(\bxi^{n})-1]
N^{\alpha}_k(\bxi^{0}) + m_k^{\alpha+1}(\bxi^{0}) N^{\alpha}_k(\bxi^{0})} 
\Delta_k(\bxi^0,\dots,\bxi^n) \\
& \overset{\eref{Def:occind}}{=} & \prod_{k=L^-}^{L^+} q^{\sum_{\alpha=1}^{n-1} 
\sum_{\beta=\alpha+1}^{n} n_k^{\beta}(\bxi^{0}) N^{\alpha}_k(\bxi^{0})} 
\prod_{\alpha=1}^n q^{-m_k^\alpha(\bxi^{n})
N^{\alpha}_k(\bxi^{0})} 
\Delta_k(\bxi^0,\dots,\bxi^n) \\
& \overset{\eref{sumrule2}}{=} & \prod_{k=L^-}^{L^+} q^{\sum_{\alpha=1}^{n-1} 
\sum_{\beta=\alpha+1}^{n} n_k^{\beta}(\bxi^{0}) N^{\alpha}_k(\bxi^{0})} 
\prod_{\alpha=1}^n q^{n^{\alpha}_k(\bxi^{0})M_k^\alpha(\bxi^{n})} 
\Delta_k(\bxi^0,\dots,\bxi^n) .
\eea

Now we compute the product over the lattice in $\bar{C}(q;\bzeta,\bfeta)$.
By sum rule \eref{sumrule3} and definition \eref{Ered} one has
\be 
\sum_{k=L^-}^{L^+} \sum_{\alpha=1}^{n-1} 
\sum_{\beta=\alpha+1}^{n} n_k^{\beta}(\bxi^{0}) N^{\alpha}_k(\bxi^{0})
= - \bar{E}(\bxi^{0})
\ee
With $\bxi^{0}=\bzeta$ and $\bxi^{n}=\bfeta$ we thus obtain \eref{aux7}. \qed

\noindent $\bullet$ {\it Step 3:} (Product expansion of symmetry generators)
Now we are in a position to prove \eref{aux5}
using the expansion \eref{dualityexpansion} for $q\neq 1$.
We find from Corollary \eref{auxn}
\be 
\bra{\bxi^{\alpha-1}} \Upsilon^+_\alpha \ket{\bxi^\alpha} =
\prod_{k=L^-}^{L^+} \bar{A}^\alpha_k(q;\bxi^{\alpha-1},\bxi^{\alpha})
\Delta^\alpha_k(\bxi^{\alpha-1},\bxi^{\alpha})
\bar{B}^\alpha_k(\bxi^{\alpha-1},\bxi^{\alpha})
%\bra{\bxi^{\alpha-1}} \rme^{S^{\alpha,+}} \ket{\bxi^{\alpha}} 
\ee
and therefore
\be 
\prod_{\alpha=1}^n \bra{\bxi^{\alpha-1}} \Upsilon^+_\alpha \ket{\bxi^\alpha} =
\prod_{k=L^-}^{L^+} 
\bar{A}_k(q;\bxi^0,\dots,\bxi^n)
\bar{B}_k(\bxi^0,\dots,\bxi^n)
\Delta_k(\bxi^0,\dots,\bxi^n)
%\bra{\bxi^{\alpha-1}} \rme^{S^{\alpha,+}} \ket{\bxi^{\alpha}} 
\ee
Observing the projector property $(\Delta_k(\bxi^0,\dots,\bxi^n))^2=
\Delta_k(\bxi^0,\dots,\bxi^n)$ we can use Lemma \eref{techlem}
and write with $\bxi^0=\bzeta$ and $\bxi^n=\bfeta$
\bea 
\bra{\bzeta} \Upsilon^+ \ket{\bfeta} 
& = & \sum_{\bxi^1  \dots \, \bxi^{n-1}}
\bar{F}(q;\bxi^0,\dots,\bxi^n) \prod_{k=L^-}^{L^+} \bar{B}_k(\bxi^0,\dots,\bxi^n)
\Delta_k(\bxi^0,\dots,\bxi^n) \nonumber \\
& \overset{\eref{aux6}}{=} &
\bar{C}(q;\bzeta,\bfeta)
\sum_{\bxi^1} \dots \sum_{\bxi^{n-1}}
\prod_{k=L^-}^{L^+} \bar{B}_k(\bxi^0,\dots,\bxi^n)
\Delta_k(\bxi^0,\dots,\bxi^n) \nonumber \\
& \overset{(\ref{aux7},\ref{theo2q1})}{=} &
%= \bar{\pi}(\bzeta)
%\prod_{k=L^-}^{L^+} \prod_{\alpha=1}^n q^{n^{\alpha}_k(\bzeta)M_k^\alpha(\bfeta)} 
%\bar{D}_1(\bzeta,\bfeta)
\bar{\pi}(\bzeta) \prod_{k=L^-}^{L^+} \prod_{\alpha=1}^n 
\left(q^{M_k^\alpha(\bfeta)}\right)^{n^{\alpha}_k(\bzeta)}
\bar{D}_{sym}(\bzeta,\bfeta)
\eea
where $\bar{D}_{sym}(\bzeta,\bfeta)$ is the reduced duality function for $q=1$.
This proves \eref{aux5} and hence
concludes the proof of Theorem \eref{Theo:duality}. \qed

\subsection{Proof of Theorem \eref{Theo:shock}}

\subsubsection{Preliminary remarks}

In \cite{Beli02} we proved for the standard ASEP $n=1$ a statement analogous to
Theorem \eref{Theo:shock} for a certain family of shock measures for the process
defined on $\Z$. The proof of \cite{Beli02} consists in three steps: (i) One 
proves, using the quantum
algebra symmetry,  that the shock measures at time $t$ defined on 
a {\it finite} lattice with $L^- = -\ceil{L/2}+1$ and $L^+=\floor{L/2}$ satisfy 
for all configurations $\bfx$ whose coordinates $x_i$ exclude the
boundary sites $L^-$ and $L^+$ a 
linear evolution equation which, in the notation of the present work, can be 
written in the form
\bel{shockASEP}
\ddt \ket{\mu^L_{\bfx}(t)} = - \sum_{\bfy} G_{\bfy\bfx} \ket{\mu^L_{\bfx}(t)} +
(\hat{B} - b) \ket{\mu^L_{\bfx}(t)}.
\ee
where $\mu^L_{\bfx}$ is a shock measure for the finite
system with $L$ sites, $G_{\bfy\bfx}$ are the transition rates of an
associated shock exclusion process that has the same non-zero rates
as the original ASEP but particle-dependent hopping rates and 
inverse hopping ratio $q$, 
boundary term $\hat{B} = w(q-q^{-1}) (\hat{n}_{L^+} - \hat{n}_{L^-})$ 
and constant $b=\lim_{L\to\infty} \E^{\mu^L_{\bfx}} B$. 
(ii) Then one uses a convergence argument based on the 
coupling and the convergence theorems in \cite{Ligg85} to show that the
contribution of the boundary term to the evolution of the sequence of shock measures
$\mu^L_{\bfx}(t)$ vanishes, at fixed $t$, in the thermodynamic limit $L\to\infty$. 
(iii) Finally, standard arguments from the theory of linear ordinary differential
equations 
allow for integration of \eref{shockASEP} to yield the time evolution of the
shock measure in infinite volume according to Theorem 2 of \cite{Beli02}
which is analogous to Theorem \eref{Theo:shock}.

Steps (ii) and (iii) employ standard tools that do not rely on the specific form 
of the renormalized hopping rates and which are independent of $n$. Hence they 
can be adapted straightforwardly to the present case and are therefore 
not repeated here. It remains only to prove (i), which, generally stated, 
follows from the following ingredients: 
(a) A similarity transformation $U^n$ that relates
the generator to itself plus some boundary term $B^n$ 
(Proposition \eref{Prop:Htrafo}),
(b) an expression for suitably defined shock measures $\mu^L_{\bfx}$ in
terms of a duality function $D^\star$ and the similarity transformation $U$
(Proposition \eref{Prop:shockmeas}), 
(c) a proof that $b=\lim_{L\to\infty} \E^{\mu^L_{\bfx}} B$ does not 
depend on $\bfx$ and that the matrix $G$ has positive transition rates
and conserves probability
(for those initial configurations $\bfx$ whose coordinates 
$x^\alpha_i$ exclude the boundary sites $L^-$ and $L^+$)
(Proposition \eref{Prop:G}), and
(d) a relation analogous to \eref{shockASEP} via duality
(Proposition \eref{Prop:shockdynfin}).

Items (b) - (d) are all
generic in the sense that they can be adapted quite straightforwardly
to other dualities for interacting
particle systems. In fact, the proposition that yields item (b) 
provides a ``recipe'' for the
construction of shock measures since $D^\star$ and $U^n$ defined below 
are not unique. Item (c) includes the explicit computation 
of the matrix elements $G_{\bfy\bfx}$ which are the transition rates for the
shock exclusion process \eref{Def:shockprocess}. 

\subsubsection{Auxiliary results}

We first prove some transformation properties.

\begin{lmm}
\label{Lem:hoppingtrafo}
Let $\sigma^{\beta\alpha}_k \sigma^{\alpha\beta}_{k+1}$ and $\hat{n}^\alpha_k$ 
be the matrices defined in \eref{Def:singlesite} and \eref{Def:localop} and 
let $\hat{\pi}$ be the matrix form of the reversible measure \eref{revmeas}. 
Under the transformations
\bel{VGamma}
\hat{V}^\alpha :=  \prod_{k=L^-}^{L^+} q^{\hat{N}_k^\alpha}, \quad 
\Gamma  := \prod_{k=L^-}^{L^+} \gamma_k
\ee
with the diagonal particle balance operators derived from \eref{Def:Nalphafun}
and the cyclic flip matrices $\gamma$ defined in \eref{cyclic}
one has 
\bel{hoppingtrafo}
\hat{V}^\gamma \sigma^{\beta\alpha}_k \sigma^{\alpha\beta}_{k+1} 
(\hat{V}^\gamma)^{-1} = \left\{ \ba{ll} 
\sigma^{\beta\alpha}_k \sigma^{\alpha\beta}_{k+1} & \gamma \neq \alpha,\beta \\
q^{2} \sigma^{\beta\alpha}_k \sigma^{\alpha\beta}_{k+1} & \gamma = \beta\\
q^{-2} \sigma^{\beta\alpha}_k \sigma^{\alpha\beta}_{k+1} & \gamma = \alpha 
\ea \right.
\ee
and
\be 
\label{Gammasigma}
\Gamma^{-1} \sigma^{\alpha\beta}_k \Gamma = \sigma^{\alpha+1,\beta+1}_{k}, \quad
\Gamma^{-1} \hat{n}^\alpha_k \Gamma = \hat{n}^{\alpha+1}_{k}, \quad
\Gamma^{-1} \hat{\pi} \Gamma = (\hat{V}^0)^{-2} \hat{\pi}
\ee
with the species indices $\alpha,\beta$ understood$\mod(n+1)$.
\end{lmm}

\Proof From the factorization 
\eref{tensorfactor} of the tensor product one has
\bel{sigmatrafo}
\hat{V}^\gamma \sigma^{\alpha\beta}_k
(\hat{V}^\gamma)^{-1} = \left\{ \ba{ll} 
\sigma^{\alpha\beta}_k & \gamma \neq \alpha,\beta \\
q^{(L^++L^--2k)} \sigma^{\alpha\beta}_k  & \gamma = \alpha\\
q^{-(L^++L^--2k)} \sigma^{\alpha\beta}_k  & \gamma = \beta 
\ea \right.
\ee
This yields \eref{hoppingtrafo}. From the definition \eref{cyclic} one finds 
\bea 
\gamma^{-1} \sigma^{\alpha\beta} \gamma & = & 
\left[|0)(n| + \sum_{\nu=1}^n |\nu)(\nu-1| \right] 
|\alpha)(\beta| \gamma \nonumber \\
& = & |\alpha+1)(\beta|\left[ |n)(0| + \sum_{\nu=1}^n |\nu-1)(\nu|\right] 
=|\alpha+1)(\beta+1|.
\eea
Similarly, the definition \eref{gammanlocal} yields $\hat{n}^\alpha \gamma= 
\gamma \hat{n}^{\alpha+1}$ from which the first two
equalities in \eref{Gammasigma} follow.

In order to prove the third equality we define
\be 
\hat{E}_{kl} := - \sum_{\alpha=1}^{n} \sum_{\beta=0}^{\alpha-1} \left(
\hat{n}_{k}^{\alpha} \hat{n}_{l}^{\beta} - 
\hat{n}_{l}^{\alpha} \hat{n}_{k}^{\beta}\right) .
\ee
From Lemma \eref{Gammasigma} we obtain
\be 
\Gamma^{-1} \hat{E}_{k,l}  \Gamma
=  \hat{E}_{k,l} 
+ \sum_{\alpha=0}^{n} 
\hat{n}_{k}^{\alpha} \hat{n}_{l}^{0}
- \sum_{\alpha=0}^{n} 
\hat{n}_{k}^{0} \hat{n}_{l}^{\alpha} + (k \leftrightarrow l)
= \hat{E}_{k,l} - 2(\hat{n}_{k}^{0} - \hat{n}_{l}^{0}) .
\ee
This, using \eref{Ered} and the resummation formula \eref{resum2},
yields the transformed energy \eref{Def:Energy} 
$\Gamma^{-1} \hat{E} \Gamma = \hat{E} - 2 \hat{E}^{0}$.
For the transformed measure this implies $\Gamma^{-1} \hat{\pi} \Gamma
= q^{-\hat{E} + 2\hat{E}^{0}}$ and with the factorization \eref{revmeasfac} 
the claim follows. \qed

The following Lemma for finite state space establishes a relation between 
duality functions and measures that is useful when the function $B$ appearing 
in the Lemma becomes irrelevant (in some sense) in the limit of infinite state 
space. Below we shall use it to construct the shock exclusion process.

\begin{lmm}
\label{Lem:bbfx}
Let $H$ be the generator of a process $\omega_t$ defined on a finite state 
space $\Omega$ which is self-dual w.r.t. a duality function $D(\xi,\omega)$
and which satisfies the intertwining relation $UD^T H = (H+B)UD^T$
for a pair of matrices $U, B$ 
 such that 
the duality matrix
$D$ and the transformation matrix $U$ satisfy
$\Phi(\omega) := \bra{s} U D^T \ket{\omega} > 0$ for all $\omega \in \Omega$.
Define the family of measures $\mu^L_{\omega}(\xi)
:= (\Phi(\omega))^{-1} \bra{\xi} U D^T \ket{\omega}$ indexed by $\omega$ 
and the functions
\bel{bt}
b(\omega) := \bra{s} \hat{\Phi} H \hat{\Phi}^{-1} \ket{\omega}, \quad
B(\omega) := \sum_{\xi\in\Omega} B_{\xi,\omega}.
\ee
One has
\bel{bt2}
U H^T = (H+B) U, \quad \E^{\mu^L_{\omega}} B = b(\omega)
\ee
for all $\omega\in \Omega$.
\end{lmm}

\Proof
The first equality follows directly from self-duality $H^T D^T = D^T H$ and
the intertwining relation.
A duality matrix can always be written in the form $D = \sum_{\omega\in\Omega} 
\ket{\omega}\bra{s} \hat{D}_\omega$ with diagonal matrix $\hat{D}_\omega$ 
that has diagonal elements $(\hat{D}_\omega)_{\xi\xi}=D(\xi,\omega)$. Also, 
$(\Phi(\omega))^{-1} > 0$ and 
$\sum_{\xi\in\Omega}\mu^L_{\omega}(\xi)=1$ so that the measures 
$\mu^L_{\omega}$ are well-defined for all $\omega \in \Omega$.
Observe that $\bra{s} B = \bra{s} \hat{B}$ for the diagonal matrix representation
\eref{diagf} of the function $B(\omega)$.
The representation
\eref{unitmatrix} of the unit matrix and conservation of probability
$\bra{s} H=0$ then leads to the chain of equations
\bea
b(\omega) & = & \Phi^{-1}(\omega) \bra{s} \hat{\Phi} H \ket{\omega} \nonumber \\
& = & \Phi^{-1}(\omega)\sum_{\xi} \Phi(\xi) \bra{\xi} H \ket{\omega}\nonumber \\
& = & \Phi^{-1}(\omega)\sum_{\xi} \bra{s} U D^T \ket{\xi} \bra{\xi} H \ket{\omega}\nonumber \\
& = &  \Phi^{-1}(\omega)\bra{s} U D^T H \ket{\omega}\nonumber \\
& = &  \Phi^{-1}(\omega)\bra{s} (H+B) U D^T \ket{\omega}\nonumber \\
& = &  \Phi^{-1}(\omega)\bra{s} \hat{B} U  D^T \ket{\omega}\nonumber \\
& = &  \bra{s} \hat{B} \ket{\mu^L_{\omega}} = \E^{\mu^L_{\omega}} B
\eea
which is the assertion of the Lemma. \qed

\subsubsection{Main part of the proof}

Items (a) - (d) form the following main part of the proof.

\begin{prop} 
\label{Prop:Htrafo}
(Item (a)) For $\gamma \in \S_{0,n}$ define the diagonal boundary matrix
\bel{boundarymatrix}
\hat{B}^\gamma := w  (q - q^{-1}) (\hat{n}^{\gamma}_{L^+} - 
\hat{n}^{\gamma}_{L^-}  ).
\ee
and let $\hat{V}^\alpha$ and $\Gamma$ be the transformation matrices \eref{VGamma}.
Then under the composite transformation
\bel{Un}
U^n = \hat{\pi} \hat{V}^n \Gamma
\ee
the generator $H$ \eref{H} on the $n$-species priority ASEP 
satisfies
\be 
\label{result2at}
H^T 
= (U^n)^{-1} \left(
H + \hat{B}^n \right) U^n.
\ee
\end{prop}

\Proof
The matrices $\hat{\pi}$, $\hat{V}^n$ and $\hat{B}^n$ are all diagonal. Therefore
$(U^n)^{-1} \hat{B}^n U^n = \Gamma^{-1}\hat{B}^n \Gamma$. Then time reversal 
via the diagonal
matrix form $\hat{\pi}$ of the reversible measure \eref{revmeas} and
transposition reduces
\eref{result2at} to 
\bel{resultv3}
H = \Gamma^{-1} \left(
\hat{V}^n H (\hat{V}^n)^{-1} + \hat{B}^n \right)  \Gamma .
\ee
which we now prove with the help of the
the decompositions
$H = H^{f} + H^{(n-1)} = H^{0} + \bar{H}$
with the corresponding bond hopping matrices 
\bel{hfk}
h^f_{k,k+1} = - w \sum_{\beta=0}^{n-1} \left[
q (\sigma^{\beta n}_k \sigma^{n \beta}_{k+1} 
- \hat{n}^{n}_k \hat{n}^{\beta}_{k+1}) 
+ q^{-1} (\sigma^{n \beta}_k \sigma^{\beta n}_{k+1} 
- \hat{n}^{\beta}_k \hat{n}^{n}_{k+1}) \right]
\ee
for the first-class particles (species $n$) and
\bel{hn1k}
h^{(n-1)}_{k,k+1} = - w \sum_{\alpha=1}^{n-1} \sum_{\beta=0}^{\alpha-1} \left[
q (\sigma^{\beta \alpha}_k \sigma^{\alpha \beta}_{k+1} 
- \hat{n}^{\alpha}_k \hat{n}^{\beta}_{k+1}) 
+ q^{-1} (\sigma^{\alpha \beta}_k \sigma^{\beta \alpha}_{k+1} 
- \hat{n}^{\beta}_k \hat{n}^{\alpha}_{k+1}) \right]
\ee
for the species of lower class and similarly
\bel{h0k}
h^0_{k,k+1} = - w \sum_{\alpha=1}^{n} \left[
q (\sigma^{0 \alpha}_k \sigma^{\alpha0}_{k+1} 
- \hat{n}^{\alpha}_k \hat{n}^{0}_{k+1}) 
+ q^{-1} (\sigma^{\alpha 0}_k \sigma^{0 \alpha}_{k+1} 
- \hat{n}^{0}_k \hat{n}^{\alpha}_{k+1}) \right]
\ee
for the vacancies (species $0$) and
\bel{hbark}
\bar{h}_{k,k+1} = - w \sum_{\alpha=2}^{n} \sum_{\beta=1}^{\alpha-1} \left[
q (\sigma^{\beta \alpha}_k \sigma^{\alpha \beta}_{k+1} 
- \hat{n}^{\alpha}_k \hat{n}^{\beta}_{k+1}) 
+ q^{-1} (\sigma^{\alpha \beta}_k \sigma^{\beta \alpha}_{k+1} 
- \hat{n}^{\beta}_k \hat{n}^{\alpha}_{k+1}) \right].
\ee

From Lemma \eref{Lem:hoppingtrafo} one has 
$\hat{V}^n h^{(n-1)}_{k,k+1}(q) (\hat{V}^n)^{-1} = h^{(n-1)}_{k,k+1}(q)$ 
and $\Gamma^{-1} h^{(n-1)}_{k,k+1}(q) \Gamma = \bar{h}_{k,k+1}(q)$ and
therefore $\Gamma^{-1} 
\hat{V}^n H^{(n-1)} (\hat{V}^n)^{-1} \Gamma = \bar{H}$.

Now we compute the transformation of the first-class part. 
In the following computations we write the 
$q$-dependence of the generator explicitly and define the transformed generator
$\tilde{H}^f(q) := \hat{V}^n H^f(q) (\hat{V}^n)^{-1}$.
Using again Lemma \eref{Lem:hoppingtrafo} yields
\be 
\tilde{h}^f_{k,k+1}(q) = - w \sum_{\beta=0}^{n-1} \left[
q (\sigma^{n \beta}_k \sigma^{\beta n}_{k+1}   
- \hat{n}^{n}_k \hat{n}^{\beta}_{k+1}) 
+ q^{-1} ( \sigma^{\beta n}_k \sigma^{n \beta}_{k+1} 
- \hat{n}^{\beta}_k \hat{n}^{n}_{k+1}) \right]
\ee
which is not stochastic. Observe, however, that we can write
$\tilde{h}^f_{k,k+1}(q) = h^f_{k,k+1}(q^{-1}) -w (q-q^{-1}) (\hat{n}^{n}_{k+1}
- \hat{n}^{n}_k )$.
Using the telescopic property of the sum this yields
$\tilde{H}^{f}(q) = H^{f}(q^{-1}) - \hat{B}^n$.
Next we apply the transformation $\Gamma$. Lemma \eref{Lem:hoppingtrafo}
yields $\Gamma^{-1} H^{f}(q^{-1}) \Gamma = H^0(q)$ and therefore
$\Gamma^{-1} 
\hat{V}^n H^{f} (\hat{V}^n)^{-1} \Gamma = H^0 - \Gamma^{-1} \hat{B}^n \Gamma$
which proves \eref{resultv3}.
\qed

Next we express the shock measure in terms of the duality function and the
similarity transformation $U^n$.

\begin{prop}
\label{Prop:shockmeas}
(Item (b)) 
Let $D^\star(\bfx,\bfeta)= 
\prod_{\alpha=2}^n \delta_{N^\alpha(\bfeta),N^\alpha(\bfx)} 
\,q^{\lambda N^0(\bfeta)} D(\bfx,\bfeta)$ for 
$\bfx \in V^L_N, \, \bfeta \in \S_{0,n}^L$ be the duality function
with $D(\bfx,\bfeta)$ given by \eref{D}.
Then with the transformation \eref{Un} and the normalization constant
\bel{Phix}
\Phi(\bfx) := \bra{s} U^n (D^\star)^T \ket{\bfx}
\ee
the shockmeasure \eref{shockmeas} can be written as
\bel{shockmeasduality}
\ket{\mu^L_{\bfx}} = \Phi(\bfx)^{-1}  U^n (D^\star)^T \ket{\bfx}.
\ee
for any $L$.
\end{prop}

\Proof
It is convenient to express the transformation \eref{Un} in the alternative
form 
\bel{Unalt}
U^n = \Gamma \hat{\pi} (\hat{V}^0)^{-1} = \Gamma \hat{\bar{\pi}}
\ee
where the first equality follows from Lemma \eref{Lem:hoppingtrafo} and the second
from the definition \eref{Ered} and the factorization \eref{revmeasfac}.
The duality matrix reads
%\bel{dualitymatrix2}
$D^\star = \sum_{\bfx}  \ket{\bfx}\bra{s} 
\hat{\tilde{Q}}_{\bfx} \hat{\wp}^2_\bfx  q^{\lambda \hat{N}^0}$
%\ee
which gives 
$(D^\star)^T \ket{\bfx} =
\hat{Q}^{\vec{c}}_{\bfx} \hat{\wp}^2_\bfx q^{\lambda \hat{N}^0}  \ket{s}$.

Now we define the subsets
$\omega^\alpha_\bfx := \underset{\beta\in\{\alpha,\dots,n\}}{\cup} 
\{\vec{x}^\beta\} \subset \Lambda^L$ 
of all particle coordinates of species $\beta \geq \alpha$ and
their complements 
$\overline{\omega^\alpha_\bfx} := \Lambda^L \setminus \omega^\alpha_\bfx$.
We also define the (unnormalized) product measures
\be  
\ket{s^{0,1}_{\bfx^{\alpha,n}}}  := 
\prod_{k\in \overline{\omega^\alpha_\bfx}} (\hat{n}^0_k + \hat{n}^1_k)
\prod_{\beta=\alpha}^{n}  \prod_{i=1}^{N^\beta(\bfx)} 
\hat{n}^\beta_{x^\beta_i} \ket{s}
\ee
with the conventions
%\be 
$\ket{s^{0,1}} := \prod_{k\in \Lambda} (\hat{n}^0_k + \hat{n}^1_k) \ket{s}$,
$\ket{s^{0,1}_{\vec{x}}} := \ket{s^{0,1}_{\bfx^{1,n}}}$
%\ee
and note that
\bel{s01x} 
\ket{s^{0,1}_{\bfx^{\beta,n}}} = \prod_{\beta=\alpha}^{n}  
\prod_{i=1}^{N^\beta(\bfx)} \sigma^{\beta 1} _{x^\beta_i} \ket{s^{0,1}}
\ee
with the raising operators $\sigma^{\beta 1}$ defined in \eref{cyclic}.
Dividing by a normalization factor $(1/2)^{|\overline{\omega^\alpha_\bfx}|}$
these are product measures with marginals $\delta_{\eta_k,\eta_k(\bfx)}$
for the sites $k\in \omega^\alpha_\bfx \subset \Lambda^L$ occupied by 
particles of species 
$\beta \geq \alpha$
and marginals $(\delta_{\eta_k,0}+\delta_{\eta_k,1})/2$ for the remaining 
sites $\overline{\omega^\alpha_\bfx}\subset \Lambda^L$.

The first observation is that 
\be 
\prod_{\alpha=2}^{n}  \prod_{i=1}^{N^\alpha(\bfx)} \hat{m}^\alpha_{x^\alpha_i}
\hat{\wp}^2_\bfx \ket{s} = \ket{s^{0,1}_{\bfx^{2,n}}}
\ee
which is a consequence of the projector property of 
$\hat{m}^\alpha_{x^\alpha_i}$ and the projection on $N^\alpha(\bfx)$ particles
for $2\leq\alpha\leq n$.
Therefore with 
$(D^\star)^T \ket{\bfx} = \hat{Q}^{\vec{c}}_{\bfx} 
\hat{\wp}^2_\bfx q^{\lambda \hat{N}^0} \ket{s}$ we have
\bea 
(D^\star)^T \ket{\bfx} 
& = &
q^{\lambda \hat{N}^0} \prod_{i=1}^{N^1(\bfx)} 
q^{-(1+c_1)\sum_{l=L^-}^{x^1_i-1}(1-\hat{m}_l^1) + 
(1-c_1)\sum_{l=x^1_i+1}^{L^+}(1-\hat{m}_l^1)} \hat{m}^1_{x^1_i} \nonumber \\
& & \times \prod_{\alpha=2}^{n} \prod_{i=1}^{N^\alpha(\bfx)}
q^{-(1+c_\alpha)\sum_{l=L^-}^{x^\alpha_i-1}(1-\hat{m}_l^\alpha) + 
(1-c_\alpha)\sum_{l=x^\alpha_i+1}^{L^+}(1-\hat{m}_l^\alpha)} 
\ket{s^{0,1}_{\bfx^{2,n}}} \nonumber \\
& = &
q^{\lambda \hat{N}^0}\prod_{i=1}^{N^1(\bfx)} 
q^{-(1+c_1)\sum_{l=L^-}^{x^1_i-1}\hat{n}^0_l + 
(1-c_1)\sum_{l=x^1_i+1}^{L^+}\hat{n}^0_l} \hat{n}^1_{x^1_i} \nonumber \\
& & \times \prod_{\alpha=2}^{n} \prod_{i=1}^{N^\alpha(\bfx)}
q^{-(1+c_\alpha)\sum_{l=L^-}^{x^\alpha_i-1}(\hat{n}_l^0+\hat{n}_l^1) + 
(1-c_\alpha)\sum_{l=x^\alpha_i+1}^{L^+}(\hat{n}_l^0+\hat{n}_l^1)} 
\ket{s^{0,1}_{\bfx^{2,n}}} \nonumber \\
& = &
q^{\lambda \hat{N}^0}\prod_{\alpha=1}^{n} \prod_{i=1}^{N^\alpha(\bfx)} 
q^{-(1+c_\alpha)\sum_{l=L^-}^{x^\alpha_i-1}\hat{n}^0_l + 
(1-c_\alpha)\sum_{l=x^\alpha_i+1}^{L^+}\hat{n}^0_l} \nonumber \\
& & \times \prod_{\alpha=2}^{n} \prod_{i=1}^{N^\alpha(\bfx)}
q^{-(1+c_\alpha)\sum_{l=L^-}^{x^\alpha_i-1}\hat{n}_l^1 + 
(1-c_\alpha)\sum_{l=x^\alpha_i+1}^{L^+}\hat{n}_l^1} 
\prod_{i=1}^{N^1(\bfx)} \hat{n}^1_{x^1_i} 
\ket{s^{0,1}_{\bfx^{2,n}}} \nonumber \\
& = &
q^{\lambda \hat{N}^0}\prod_{\alpha=1}^{n} \prod_{i=1}^{N^\alpha(\bfx)} 
q^{-(1+c_\alpha)\sum_{l=L^-}^{x^\alpha_i-1}\hat{n}^0_l + 
(1-c_\alpha)\sum_{l=x^\alpha_i+1}^{L^+}\hat{n}^0_l} \nonumber\\
& & \times \prod_{k\in\omega^2_\bfx} q^{-\hat{N}^1_k} \prod_{\alpha=2}^{n}
q^{-c_\alpha N^\alpha(\bfx) \hat{N^1}}
%\prod_{i=1}^{N^1(\bfx)} \hat{n}^1_{x^1_i} \ket{s^{0,1}_{\bfx^{2,n}}}.
\ket{s^{0,1}_{\vec{x}}}.
\eea
The substitution of $\hat{m}^1$ by $\hat{n}^1$ in the second equality comes 
from the fact that $\{\vec{x}^1\}$ is in the complement
$\overline{\omega^2_\bfx}$ where one has $\hat{n}^\alpha_y \ket{s^{0,1}_\bfx} = 0$ 
for $2 \leq \alpha \leq n$. In the last equality we have used that by construction
$\prod_{i=1}^{N^1(\bfx)} \hat{n}^1_{x^1_i} \ket{s^{0,1}_{\bfx^{2,n}}}=
\ket{s^{0,1}_{\vec{x}}}$.

The next step is to compute $U^n (D^\star)^T \ket{\bfx}$
using \eref{Unalt}. To this end
we make a decomposition of $\bar{E}(\cdot)$ \eref{Ered} as follows. 
From \eref{Def:Eab} we construct
the doubly reduced energy
\be 
\bar{\bar{E}}(\cdot) := 
\sum_{\alpha=3}^{n} \sum_{\beta=2}^{\alpha-1} E^{\alpha\beta}(\cdot).
\ee
The representation \eref{Eab} of the partial energies allows us to write
\bel{Ered2}
\bar{E}(\cdot) = \bar{\bar{E}}(\cdot) + \sum_{\alpha=2}^{n} E^{\alpha 1}(\cdot)
= \bar{\bar{E}}(\cdot) -  \sum_{k=L^-}^{L^+} m_k^2(\cdot) N_k^1(\cdot)
\ee
with the corresponding decomposition of the matrix form of the 
reduced measure \eref{pired}
\bel{pired3}
\hat{\bar{\pi}} = \hat{\bar{\bar{\pi}}} q^{\sum_{k=L^-}^{L^+} 
\hat{m}_k^2 \hat{N}_k^1}.
\ee
Since $\hat{\bar{\bar{\pi}}}\ket{s^{0,1}_{\bfx^{2,n}}} 
= \bar{\bar{\pi}}(\bfx) \ket{s^{0,1}_{\bfx^{2,n}}}$ and likewise 
$\hat{\bar{\bar{\pi}}}\ket{s^{0,1}_{\vec{x}}} = 
\bar{\bar{\pi}}(\bfx) \ket{s^{0,1}_{\vec{x}}}$
we find
\be 
\hat{\mathcal{P}}_\bfx \hat{\bar{\pi}} \ket{s^{0,1}_{\vec{x}}} 
= \bar{\bar{\pi}}(\bfx)
\prod_{k\in\omega^2_\bfx} q^{\hat{N}_k^1} \ket{s^{0,1}_{\vec{x}}}
\ee
which comes from the fact that $m_k^2(\bfx)=\sum_{x\in\omega^2_\bfx} \delta_{k,x}$.
Therefore we arrive at the intermediate result
\bea
\hat{\bar{\pi}} (D^\star)^T \ket{\bfx} 
& = & \bar{\bar{\pi}}(\bfx) 
q^{\lambda \hat{N}^0}\prod_{\alpha=1}^{n} \prod_{i=1}^{N^\alpha(\bfx)} 
q^{-(1+c_\alpha)\sum_{l=L^-}^{x^\alpha_i-1}\hat{n}^0_l + 
(1-c_\alpha)\sum_{l=x^\alpha_i+1}^{L^+}\hat{n}^0_l} \nonumber \\
& & \times  \prod_{\alpha=2}^{n}
q^{-c_\alpha N^\alpha(\bfx) \hat{N^1}}
%\prod_{i=1}^{N^1(\bfx)} \hat{n}^1_{x^1_i} \ket{s^{0,1}_{\bfx^{2,n}}}.
\ket{s^{0,1}_{\vec{x}}}.
\eea

Now we consider the transformation $\Gamma$. With the definition
%\be 
$\ket{s^{n,0}} := \prod_{k=L^-}^{L^+} (\hat{n}^0_k+\hat{n}^n_k) \ket{s}$
%\ee
we obtain from \eref{s01x} 
\be 
\ket{s^{n,0}_{\bfx^{-}}} := \Gamma \ket{s^{0,1}_\bfx} =
\prod_{i=1}^{N^1(\bfx)} \hat{n}^0_{x^1_i} 
\prod_{\alpha=2}^{n}  \prod_{i=1}^{N^{\alpha}(\bfx)} 
\sigma^{\alpha-1,+}_{x^{\alpha}_i} \ket{s^{n,0}}
\ee
with the representation matrices $\sigma^{\alpha,+}$ \eref{fundgln}.
Notice that this is a (unnormalized) product measure with particles
of type $\alpha-1$ at the positions $x^{\alpha}_i$ with probability 1, 
corresponding to shock markers of type $\alpha$ at these positions.
Putting these results together yields
\bea 
U^n (D^\star)^T \ket{\bfx} 
& = & \bar{\bar{\pi}}(\bfx)
\Gamma q^{\lambda \hat{N}^0} \prod_{\alpha=1}^{n} \prod_{i=1}^{N^\alpha(\bfx)} 
q^{-(1+c_\alpha)\sum_{l=L^-}^{x^\alpha_i-1}\hat{n}^0_l + 
(1-c_\alpha)\sum_{l=x^\alpha_i+1}^{L^+}\hat{n}^0_l} \Gamma^{-1} \nonumber \\
& & \times \Gamma \prod_{\alpha=2}^{n}
q^{-c_\alpha N^\alpha(\bfx) \hat{N^1}}
\Gamma^{-1} \ket{s^{n,0}_{\bfx^{-1}}} \nonumber \\
& = & \bar{\bar{\pi}}(\bfx) q^{\lambda \hat{N}^n} 
\prod_{\alpha=1}^{n} \prod_{i=1}^{N^\alpha(\bfx)} 
q^{-(1+c_\alpha)\sum_{l=L^-}^{x^\alpha_i-1}\hat{n}^n_l + 
(1-c_\alpha)\sum_{l=x^\alpha_i+1}^{L^+}\hat{n}^n_l}  \nonumber \\
& & \times \prod_{\alpha=2}^{n}
q^{-c_\alpha N^\alpha(\bfx) \hat{N^0}}
\ket{s^{n,0}_{\bfx^{-}}}.
\eea

Finally we set $c_\alpha = 0$ for all $\alpha$ which leads to
\bea 
U^n (D^\star)^T \ket{\bfx} 
& = & \bar{\bar{\pi}}(\bfx) q^{\lambda \hat{N}^n} 
\prod_{k=1}^{N(\bfx)}  q^{-\hat{N}^n_k} \ket{s^{n,0}_{\bfx^{-}}}.
\eea
and thus proves that
$U^n (D^\star)^T \ket{\bfx}/\Phi(\bfx)$ with $\Phi(\bfx) = \bra{s} U^n (D^\star)^T 
\ket{\bfx}$ is a product measure with shock markers of type 
$\alpha$ at positions $x^{\alpha}_i$ (which are particles of species 
$\alpha-1$).

In order to identify this measure with the shock measure \eref{shockmeas}
we compute the normalization $\Phi(\bfx)$. 
Fixing $K=N(\bfx)$ we have
\bea 
\prod_{k=1}^{N(\bfx)}  q^{-\hat{N}^n_k}
& = & \prod_{i=1}^K 
q^{-\sum_{l=L^-}^{x_i-1}\hat{n}^n_l  + 
\sum_{l=x_i+1}^{L^+}\hat{n}^n_l} \nonumber \\
& = & \prod_{k=L^-}^{x_1-1} q^{-K\hat{n}^n_k}
\prod_{k=x_1+1}^{x_2-1} q^{[-(K-1)+1]\hat{n}^n_k}  \dots 
\prod_{k=x_{i}+1}^{x_{j+1}-1} q^{(2i-K)\hat{n}^n_k} \dots
\prod_{k=x_{K}+1}^{L^+} q^{K\hat{n}^n_k} \nonumber \\
& = &
\prod_{i=0}^{K} \prod_{k=x_{i}+1}^{x_{i+1}-1} q^{(2i-K)\hat{n}^n_k}
\eea
with the conventions $x_0=L^--1$ and $x_{K+1} = L^++1$.

Using the product structure of $\ket{s^{n,0}_{\bfx^{-}}}$ one thus gets
the normalization
\bel{Phi} 
\Phi(\bfx) = \bar{\bar{\pi}}(\bfx)
\prod_{i=0}^{K} (1+q^{2i-K+\lambda})^{x_{i+1}-x_i-1}.
\ee
and the marginal densities
\bel{shockmargphi}
\rho_i = \frac{q^{2i-K+\lambda}}{1+q^{2i-K+\lambda}}
\ee
for the sites $x_i < k < x_{i+1}$ between the shock markers $i$ and
$i+1$ which are those defined in \eref{shockmargfunction}. \qed

The third ingredient (c) establishes the link between the matrix elements
$\hat{\Phi} H \hat{\Phi}^{-1}$ and the shock exclusion process.

\begin{prop}
\label{Prop:G}
(Item (c)) Let $\tilde{V}^L_N := \{\bfx \in V^L_N \ | \ L^- < x^\alpha_i < L^+ \ 
\forall \ x^\alpha_i \in \{\bfx\} \}$ 
be the set of configurations with all particle positions restricted to the segment
$[L^-+1,L^+-1]$ of $\Lambda_L$ and define for $\bfx \in V^L_N$ the function
$b(\bfx):=\E^{\mu^L_{\bfx}} B^n$ for the shock measures 
$\mu^L_{\bfx}$ \eref{shockmeas} and boundary matrix \eref{boundarymatrix}.
The following holds:\\

\noindent (a) $b(\bfx) = b$ for all $\bfx \in \tilde{V}^L_N$ with 
$b = w(q-q^{-1})(\rho^+-\rho^-)$ and $\rho^\pm = \E^{\mu^L_{\bfx}} n_{L^\pm}$.\\

\noindent (b) Define the matrix
\bel{shockdyn2}
G := \hat{\Phi} H \hat{\Phi}^{-1} - b
\ee
with the shock
normalization \eref{Phi}. Then for the negative off-diagonal matrix elements
$w_{\bfy,\bfx} := -G_{\bfy\bfx}$ one has positivity $w_{\bfy,\bfx}\geq 0$
for all $\bfx,\bfy \in V^L_N$, $\bfy\neq\bfx$ 
and conservation of probability 
$\sum_{\bfy\in V^L_N} G_{\bfy\bfx} =0$ for $\bfx \in \tilde{V}^L_N$.\\

\noindent (c) The negative off-diagonal matrix elements
$w_{\bfy,\bfx}$ are the the shock transition 
rates \eref{shockhoppingrates}.
\end{prop}

\Proof
Part (a) is trivial for the product measure \eref{shockmeas} 
since for any $\bfx$ with coordinates $x_i^\alpha \in [L^-+1,L^+-1]$
its boundary marginals and hence
the expectation $\E^{\mu^L_{\bfx}} B^\gamma$ does not depend on $\bfx$
for any $\gamma \in \S_{0,n}$. The value
of $b$ follows directly from the definition of the boundary matrix $B^n$ 
\eref{boundarymatrix}.

In order to prove part (b) we first note that positivity is trivial since by 
construction $\Phi(\bfx)>0$ and $H_{\bfy\bfx} \leq 0$ for all non-equal pairs 
$\bfx,\bfy \in V_{N}^L$. Applying Lemma \eref{Lem:bbfx} 
to the present setting then yields conservation of 
probability for configurations $\bfx \in \tilde{V}^L_N$. 

The last item (c) to be proved is the identification of
$-G_{\bfy\bfx}=-\Phi(\bfy)H_{\bfy\bfx}(\Phi(\bfx))^{-1}$ with 
the shock transition 
rates \eref{shockhoppingrates}. We recall that due to the transformation
$\Gamma$ in the definition of the shock measure a configuration
$\bfx = (\vec{x},\vec{\alpha})$ of the
$n$-species priority ASEP corresponds to a configuration of shock markers of 
type $\alpha_i-1$ at the sites $x_i$. 
Observing that $G_{\bfy\bfx} = 0$ $\iff$ $H_{\bfy\bfx} = 0$  
we need to consider only configurations $\bfy$ that differ from $\bfx$
either by a single unit displacement (i.e., $y^\alpha_l = x^\alpha_l\pm 1$ 
for some specific $x^\alpha_i$ and
$y^{\alpha'}_{l'} = x^{\alpha'}_{l'}$ for all other coordinates) 
or by an interchange of color (when $y^\beta_{k}=x^\alpha_{l}=x$ and
$y^\alpha_{l}=x^\beta_{k}=x+1$).

We define
\bel{Psi}
\Psi(\bfx) := \prod_{j=0}^{K} (1+q^{2j-K+\lambda})^{x_{j+1}-x_j-1}
\ee
which allows us to
split the normalization \eref{Phi} into two parts 
$\Phi(\bfx) = \hat{\hat{\pi}}(\bfx)\Psi(\bfx)$. Thus the matrix element
$G_{\bfy\bfx}$ becomes a product of three terms
\be 
G_{\bfy\bfx} = \frac{\hat{\hat{\pi}}(\bfy)}{\hat{\hat{\pi}}(\bfx)}
\times \frac{\Psi(\bfy)}{\Psi(\bfx)} \times H_{\bfy\bfx}.
\ee

First we consider the case of a jump
$x^\alpha_{j\pm 1} - x^\alpha_j >1$. 
Then $H_{\bfy\bfx} = q^{\pm 1}$ and $\hat{\hat{\pi}}(\bfy)/\hat{\hat{\pi}}(\bfx)=1$.
From \eref{Phi} and \eref{shockmargphi} one finds
\be 
\frac{\Psi(\bfy)}{\Psi(\bfx)} = \left(\frac{1+q^{2(j-1)-K+\lambda}}
{1+q^{2j-K+\lambda}}\right)^{\pm 1} = \left(\frac{1-\rho_j}
{1-\rho_{j-1}}\right)^{\pm 1}.
\ee
From \eref{shockmargphi} one also finds
\be
\frac{\rho_j}{(1-\rho_j)} = q^{2j-K+\lambda}.
\ee
Following \cite{Beli02} we note that \eref{shockfugacityratio} then yields
\bea 
(q-q^{-1})\rho_j(1-\rho_j) 
%& = &
%q \left(1-\frac{\rho_{l-1}(1-\rho_{l})}{\rho_{l}(1-\rho_{l-1})}\right)
%\rho_l(1-\rho_l) \nonumber \\
& = & q (\rho_{j}-\rho_{j-1})\frac{1-\rho_j}{1-\rho_{j-1}} \\
(q-q^{-1})\rho_{j-1}(1-\rho_{j-1}) 
%& = &
%q^{-1} \left(\frac{\rho_{l}(1-\rho_{l-1})}{\rho_{l-1}(1-\rho_{l})}-1\right)
%\rho_{l-1}(1-\rho_{l-1}) \nonumber \\
& = & q^{-1} (\rho_{j}-\rho_{j-1}) \frac{1-\rho_{j-1}}{1-\rho_j}.
\eea
Therefore
%\be \be 
%-G_{\bfy\bfx} = w q^{\pm 1} \left(\frac{1-\rho_j}
%{1-\rho_{j-1}}\right)^{\pm 1}
%\ee
$G_{\bfy\bfx} = - w q^{\pm 1} \left[(1-\rho_j)/(1-\rho_{j-1})\right]^{\pm 1}$
which gives the shock hopping rates
\bea 
w^+_j & = & w q \left(\frac{1-\rho_j}
{1-\rho_{j-1}}\right) = w \frac{(q-q^{-1})\rho_j(1-\rho_j)}{\rho_{j}-\rho_{j-1}}\\
w^-_j & = & w q^{-1} \frac{1-\rho_{j-1}}{1-\rho_j} 
= w \frac{(q-q^{-1})\rho_{j-1}(1-\rho_{j-1})}{\rho_{j}-\rho_{j-1}}
\eea
in agreement with the definition \eref{shockhoppingrates} of the shock 
exclusion process.

Next we consider the case of color exchange. In this case 
$\Phi(\bfy)/\Phi(\bfx)=1$ and proof reduces to calculating
$G_{\bfy\bfx}=H_{\bfy\bfx}\bar{\bar{\pi}}(\bfy)/\bar{\bar{\pi}}(\bfx)$.
For $2 \leq \alpha \leq n$ reversibility yields $G_{\bfy\bfx}= H_{\bfx\bfy}$.
For $\alpha=1$ one has $\bar{\bar{\pi}}(\bfy)/\bar{\bar{\pi}}(\bfx)=1$ and
therefore $G_{\bfy\bfx}= H_{\bfy\bfx}$. Both cases are in 
accordance with the colour exchange rates \eref{colourrates} 
of the shock exclusion process. \qed

The final building block in the proof of Theorem \eref{Theo:shock} reads:

\begin{prop}
\label{Prop:shockdynfin}
(Item (d)) Let $\ket{\mu^L_{\bfx}}$ be the shock measure defined in 
\eref{shockmeas} and
let $H^B$ be the evolution operator 
\bel{def:HL}
H^B := H + \hat{B}^n - b.
\ee
Then 
\bel{shockevolutionfin}
H^B \ket{\mu^L_{\bfx}} = \sum_{\bfy} G_{\bfy\bfx} \ket{\mu^L_{\bfx}}
\ee
with the transition matrix elements $G_{\bfy\bfx} = \bra{\bfy} G \ket{\bfx}$.
\end{prop}

\Proof
From selfduality \eref{Def:selfduality} with the duality matrix of Theorem
\eref{Theo:duality}, Proposition \eref{Prop:Htrafo}, and Proposition 
\eref{Prop:shockmeas} we have the intertwining relation
$(H+\hat{B}^n-b) U^n D^T =U^n D^T (H-b)$.
Thus with the representation \eref{unitmatrix} of the unit matrix and 
\eref{shockmeasduality} of the
shock measure
\bea
(H+\hat{B}^n-b) \ket{\mu^L_{\bfx}} 
& = & (\Phi(\bfx))^{-1} U^n D^T (H-b) \ket{\bfx} \\
& = & \sum_\bfy \frac{\Phi(\bfy)}{\Phi(\bfx)} \ket{\mu_\bfy} \bra{\bfy} (H-b) 
\ket{\bfx} \\
& = & \sum_\bfy \ket{\mu_\bfy} \bra{\bfy} 
(\hat{\Phi} H \hat{\Phi}^{-1}-b) \ket{\bfx}.
\eea
With the definition \eref{shockdyn2} of $G$
this proves Proposition \eref{Prop:shockdynfin}. \qed

As detailed in the preliminary remarks this completes the 
proof of Theorem \eref{Theo:shock}.

\section*{Acknowledgements}
GMS thanks the Institute of Mathematics and Statistics at the University of
S\~ao Paulo for kind hospitality and L.R.G. Fontes for 
stimulating comments. This work was supported by FAPESP (2015/15258-9),
CNPq (307347/2013-3) and DFG (SCHU 827/9-1).

\appendix

\section{Some conventions and useful formulas}
\label{A1}

The Kronecker-$\delta$ is defined by
\bel{Def:Kron}
\delta_{x,y} := \left\{ \ba{ll} 1 & \mbox{ if } x=y \\ 0 & \mbox{ else } \ea \right.
\ee
for $x,y$ from any set. 
For $x\in\R$ we define
\be 
\label{Def:Thetasignfun}
\Theta(x) := \left\{ \ba{rl} 1 & x > 0 \\ 0 & x \leq 0 \ea \right. , \quad
\sigma(x) := \Theta(x) - \Theta(-x) = \left\{ 
\ba{rl} 1 & x > 0 \\ 0 & x = 0 \\ -1 & x < 0. \ea \right.
\ee
Then one has for $k,l,m,n \in \Z$
\bel{sumTheta}
\sum_{m=k}^{l-1} \delta_{m,n} = \Theta(l-n) - \Theta(k-n)
\ee
which we take as the definition of a summation when the 
upper summation index is smaller than the lower summation index. In particular,
one has for any summable object $f_n$, $n\in\Z$,
\be 
\sum_{n=k}^{l-1} f_n = \left\{ \ba{cc} 
\displaystyle 0 & l=k \\
\displaystyle - \sum_{n=l}^{k-1} f_n &  l < k.
\ea \right.
\ee
This implies analogous relations for products of $f_n$ when the 
upper product index is smaller than the lower product index
through the formal identity 
$\prod_{n=k}^{l-1} = \exp{\left(\sum_{n=k}^{l-1} \ln{(f_n)} \right)}$.
For $k=l$ we define $\prod_{n=k}^{k-1} f_n := 1$ even
if $f_k=0$, consistent with the convention $0^0=1$.

We also note various sum rules that are used in several places in the proofs. 
For numbers or matrices $a_k,b_l$ one has
\bea
\label{sumrule1} 
\sum_{k=L^-}^{L^+} \sum_{l=L^-}^{k-1} (a_k-a_l) 
& = & \sum_{k=L^-}^{L^+} (2k-L^+-L^-) a_{k} \\
\label{sumrule2} 
\sum_{k=L^-}^{L^+} a_k \left(\sum_{l=L^-}^{k-1} b_l - \sum_{l=k+1}^{L^+} b_l\right)
& = & - \sum_{k=L^-}^{L^+} b_k \left(\sum_{l=L^-}^{k-1} a_l - \sum_{l=k+1}^{L^+} a_l\right) \\
\label{sumrule3}
\sum_{k=L^-}^{L^+} \sum_{l=k+1}^{L^+} a_k b_l 
& = & \sum_{k=L^-}^{L^+} \sum_{l=L^-}^{k-1} b_k  a_l .
\eea

As a simple consequence of these sum rules we have
\begin{lmm}
\label{Lemma:resum} (Resummation)
For numbers or matrices $a_k, b_l$ 
define 
\be 
A_k := \left(\sum_{l=L^-}^{k-1} a_l - \sum_{l=k+1}^{L^+} a_l\right), \quad
B_k := \left(\sum_{l=L^-}^{k-1} b_l - \sum_{l=k+1}^{L^+} b_l\right). 
\ee
Then one has 
\bel{resum1}
\sum_{k=L^-}^{L^+} A_k = \sum_{k=L^-}^{L^+} (L^++L^--2k) a_k =
\sum_{k=L^-}^{L^+} \sum_{l=L^-}^{k-1} (a_l-a_k)
\ee
and 
\bel{resum2}
\sum_{k=L^-}^{L^+} a_k B_k  = - \sum_{k=L^-}^{L^+} b_k A_k .
\ee
\end{lmm}

For $c\in\C$ and $q,\,q^{-1} \in \C \setminus 0$ we define the symmetric 
$q$-number by
\bel{Def:qnumber}
[c]_q := \frac{q^c - q^{-c}}{q-q^{-1}}.
\ee
For integers $n \in \N$ the $q$-factorial and the $q$-binomial
coefficient are defined by
\bel{Def:qfactorial}
[n]_q! := \prod_{k=1}^n [k]_q, \quad {n \choose k}_q := \frac{[n]_q!}{[k]_q![n-q]_q!}
\ee
and $[0]_q! := 1$. 
For finite-dimensional square matrices $A$ the expression $[A]_q$ is defined analogously to \eref{Def:qnumber} through the Taylor expansion of the exponential.

For two endomorphisms on some vector space represented by square 
matrices $A$ and $B$ we define the commutator symbol $\comm{A}{B}:= AB - BA$
with the matrix product $(AB)_{mn} = \sum_k A_{mk} B_{kn}$.
The Kronecker product $A\otimes B$ is defined for arbitrary
rectangular matrices as follows.

\begin{df}
\label{Def:Kronecker} (Kronecker product)
Let $A$ and $B$ be two matrices with $m_A$ ($m_B$) rows and
$n_A$ ($n_B$) columns with matrix elements $A_{ij}$, $1\leq i \leq m_A,
1\leq j \leq n_A$ and $B_{ij}$, $1\leq i \leq m_B,
1\leq j \leq n_B$ respectively. The Kronecker product
$A\otimes B$ is a $m_Am_B \times n_An_B$-matrix $C$ with matrix elements
$C_{kl}=A_{qp} B_{kl}$ for $(p-1)n_B+1 \leq l \leq p n_B$, 
$(q-1)m_B+1 \leq k \leq qm_B$ where $1\leq p \leq n_A$ and
$1\leq q \leq m_A$.
\end{df}

A matrix is called nilpotent of degree $k$ if $A^k=0$. Here $0$ represents
the matrix with all matrix elements $A_{mn}$ are equal to 0. A matrix $A$
is called a projector if $A^2=A$. We call a matrix $A$ satisfying $A^3=A$
a signed projector.

We mention the following simple projector lemma:
\begin{lmm}
\label{projlem} (Exponential of projectors)
(a) Let $P$ be a projector. Then for $c \in \C$ one has 
$c^P = 1 +(c-1)P$.\\
(b) Let $Q$ be a signed projector. 
Then for $c \in \C\setminus 0$ one has 
$c^Q = 1 +\frac{1}{2}(c-c^{-1})Q + 
\frac{1}{2}(c+c^{-1}-2)Q^2$.
\end{lmm}

This is an immediate consequence of the Taylor expansion of the exponential
and the projector property. In particular, we note that 
$\Theta^2(n) = \Theta(n)$ and $\sigma^3(n) = \sigma(n)$ 
so that Lemma \eref{projlem} can be applied to exponentials of these functions.

\end{document}